\documentclass[final,leqno]{siamltex}

\usepackage{ifthen}

\usepackage{amsmath,amsfonts,amssymb,algorithm,algorithmic,psfrag}
\usepackage{array}
\usepackage{psfrag}
\usepackage{subcaption}
\usepackage{bm}
\usepackage{color}
\usepackage{enumitem}

\usepackage{ifpdf}
\ifpdf
    \usepackage[pdftex]{graphicx}
    \usepackage[pdftex]{hyperref}
    \hypersetup{
        unicode=false,          
        pdftoolbar=true,        
        pdfmenubar=true,        
        pdffitwindow=false,     
        pdfstartview={FitH},    
        pdfnewwindow=true,      
        colorlinks=true,        
        linkcolor=red,          
        citecolor=blue,         
        filecolor=magenta,      
        urlcolor=cyan           
    }

\else
    \usepackage{graphicx}
    \usepackage{pstricks}
    
    \newcommand{\href}[2]{#2}
\fi

\definecolor{DkBlue}{rgb}{0,0,.5}

\newcommand{\dd} {\, \mathrm{d}}

\newcommand{\cA} {{\mathcal A}}

\newcommand{\R} {\mathbb R}

\newcommand{\cT} {\mathcal T}
\newcommand{\cS} {\mathcal S}
\newcommand{\cV} {\mathcal V}

\newcommand{\p}{\partial}

\newcommand{\cancel}[1]{}
\newcommand{\raw}{\rightarrow}

\newcommand{\eps}{\epsilon}

\newcommand\diam{\textnormal{diam}}
\newcommand\dist{\textnormal{dist}}
\newcommand{\vn}[1]{\left|\left|#1\right|\right|}
\newcommand{\vsn}[1]{\left|#1\right|}

\newcommand{\ds}{\displaystyle}

\newcommand{\bc}{{\textbf{c}}}

\newcommand{\bv}{\textbf{v}}

\newcommand{\bx}{\textbf{x}}

\newcommand{\origin}{{\textbf{o}}}

\newcommand{\hx}{\hat{x}}
\newcommand{\hy}{\hat{y}}

\newcommand{\hpn}[3]{\vn{#1}_{H^{#2}(#3)}}   
\newcommand{\hpsn}[3]{\vsn{#1}_{H^{#2}(#3)}} 
\newcommand{\wmpn}[4]{\vn{#1}_{W^{#2,#3}(#4)}}
\newcommand{\wmpsn}[4]{\vsn{#1}_{W^{#2,#3}(#4)}}

\newcommand{\cP}{{\mathcal P}}

\definecolor{DkGreen}{rgb}{0,.5,0}

\definecolor{DkGreen}{rgb}{0,.5,0}

\newtheorem{remark}[theorem]{Remark}
\newtheorem{cor}[theorem]{Corollary}

\begin{document}

\title{Interpolation Error Estimates for Harmonic Coordinates On Polytopes}

\author{Andrew Gillette\thanks{Department of Mathematics, University of Arizona, Tucson, Arizona, USA, {\tt agillette@math.arizona.edu}}
 \and 
Alexander Rand\thanks{CD-adapco, Austin, Texas, USA, {\tt alexander.rand@cd-adapco.com}}
}

\maketitle

\begin{abstract}
Interpolation error estimates in terms of geometric quality measures are established for harmonic coordinates on polytopes in two and three dimensions.
First we derive interpolation error estimates over convex polygons that depend on the geometric quality of the triangles in the constrained Delaunay triangulation of the polygon. 
This characterization is sharp in the sense that families of polygons with poor quality triangles in their constrained Delaunay triangulations are shown to produce large error when interpolating a basic quadratic function.
Non-convex polygons exhibit a similar limitation: large constrained Delaunay triangles caused by vertices approaching a non-adjacent edge also lead to large interpolation error.
While this relationship is generalized to convex polyhedra in three dimensions, the possibility of sliver tetrahedra in the constrained Delaunay triangulation prevent the analogous estimate from sharply reflecting the actual interpolation error.
Non-convex polyhedra are shown to be fundamentally different through an example of a family of polyhedra containing vertices which are arbitrarily close to non-adjacent faces yet the interpolation error remains bounded.
\end{abstract}

\section{Introduction}

Interpolation error estimates over triangles or higher-dimensional simplices are an essential component of the most common finite element analyses.
The simplest approach  involves restricting attention to simplices of bounded aspect ratio, which, in two-dimensions, is equivalent to the familiar minimum angle condition.
However, this is an overly restrictive approach and identical estimates have been established over broader classes of simplices.
In two-dimensions, improved estimates involve the maximum angle condition~\cite{BA76,Kr91,GMW99} which can be generalized to higher dimensions~\cite{Ja76,Kr92,Sh94}.
These types of conditions and ideas can then be applied to 
estimates in different norms~\cite{Sy57,Ac01,Ra12}, 
other types of elements~\cite{Du99,AD99} or, 
in some sense, even to finite volume methods~\cite{MTTW99}.
The interpolation estimates are sharp in the sense that there exist sequences of function/simplex pairs that realize the estimate~\cite{BA76,Ap99}.
While the maximum angle condition is not strictly necessary for convergence of finite element methods~\cite{HKK12}, 
the known counterexamples involve finite element spaces which contain subspaces corresponding to shape-regular meshes.
Thus, some form of shape regularity still appears to be a necessary ingredient for successful interpolation.

Polygonal finite elements can be established as a generalization of triangular finite elements by attempting to preserve many properties from the triangular case.
This essential framework has been applied to a variety of different finite element contexts~\cite{SMSB2001,ESW2006,TS06,WBG07,MKBWG2008}  and shares much of the underlying theory with related numerical methods for polygonal meshes: 
mimetic methods~\cite{AVZ2013,BBM2013,LMS2014,BLS2005,BH2006}, 
virtual element methods~\cite{dBCMMR2013,dVBMR2014,dBM2013,CMRS2014,MRS2014}, 
weak Galerkin methods~\cite{W2014,WY2013,MWY2015},
compatible discretization operator schemes~\cite{BE2014a,BE2014b,LCAV2014}, 
and some ``meshfree'' methods~\cite{LL02,SW2007,BBO02}.

Generalized barycentric coordinates (GBCs) are a common approach to defining basis functions for polygonal and polyhedral finite element methods.
A set of GBCs on an $n$-dimensional polytope will \emph{include} in their span the set of linear functions on $\R^n$, but the GBCs themselves are in general not polynomials.
Construction of GBCs is not unique and many competing variants have been developed.
These include
Wachspress coordinates~\cite{W1975},
mean value coordinates~\cite{F2003},
Sibson coordinates~\cite{S1980},
maximum entropy coordinates~\cite{Su04},
Poisson coordinates~\cite{LH2013},
and harmonic coordinates~\cite{JMRGS07,MKBWG2008}.
Motivation for GBCs include data interpolation~\cite{S1980}, computer graphics~\cite{JMRGS07}, and mesh generation~\cite{CT14}, as well as finite element applications.
Interpolation error, the focus of the paper, requires the development uniform estimates over some class of polytopes for a particular generalized barycentric construction. 
Since all GBCs do not admit estimates over the same geometric restrictions, identification of coordinates which admit the widest set of polytopes provides one method of differentiating particular constructions.

Among all possible GBCs, harmonic coordinates are optimal in the context of interpolation error: 
interpolation by harmonic coordinates minimizes the $H^1$-seminorm over all functions satisfying the requisite boundary data.
Thus, a bound on the interpolation error of harmonic coordinates in terms of geometric properties of the domain serves as a benchmark for the error analysis of all other types of GBCs.
Generalizing geometric conditions from simplices to polytopes, however, must be done with some care; some conditions which are equivalent on triangles are no longer equivalent on polygons.
For example, requiring a minimum angle is equivalent to bounded aspect ratio on triangles, while a generic convex polygon can have no small angles but poor aspect ratio, as in the case of a rectangle with a large ratio of width to height.
Prior error analysis is primarily aimed at generalizing the minimum angle condition on triangles. 
For harmonic coordinates error estimates were established for convex polygons with bounded aspect ratio~\cite{GRB2011}, and for other GBCs, some additional restrictions are were placed~\cite{GRB2011,RGB2011b}.
In dimension larger than two, similar estimates have been established for Wachspress coordinates for some families of shape-regular polytopes~\cite{FGS2013}.

\begin{figure}[ht]
   \begin{center}
  \begin{subfigure}[b]{0.18\textwidth}
    \includegraphics[height=\textwidth]{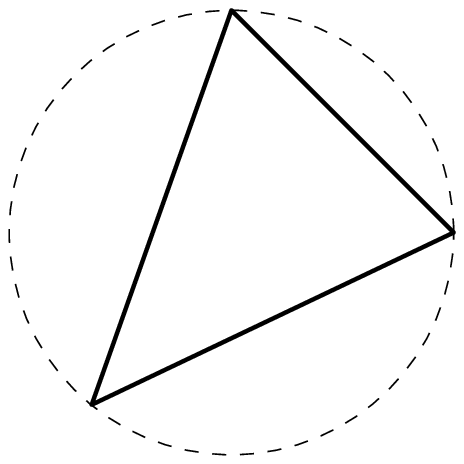}
    \caption{}
    \label{fg:triangleQualityA}
  \end{subfigure}
\hspace{0.04\textwidth}
  \begin{subfigure}[b]{0.18\textwidth}
    \includegraphics[height=\textwidth]{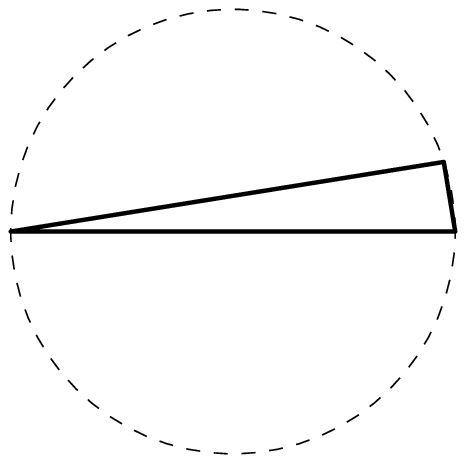}
    \caption{}
    \label{fg:triangleQualityB}
  \end{subfigure}
  \begin{subfigure}[b]{0.18\textwidth}
    \hspace{.15\textwidth}
    \includegraphics[height=\textwidth]{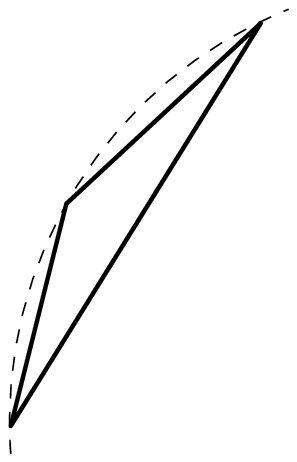}
    \caption{}
    \label{fg:triangleQualityC}
  \end{subfigure}
       \end{center}
\caption{Three types of triangles from the perspective of interpolation quality metrics: 
(A) a good aspect ratio triangle, 
(B) a poor aspect ratio triangle with a single small angle (and thus modest circumradius) and 
(C) a poor aspect ratio triangle with a large angle (and thus large circumradius).
Standard interpolation error estimates hold for the first two types but fail for the third.
}\label{fg:triangleQuality}
\end{figure}

\begin{figure}[t]
\begin{center}
  \begin{subfigure}[b]{0.18\textwidth}
    \includegraphics[height=\textwidth]{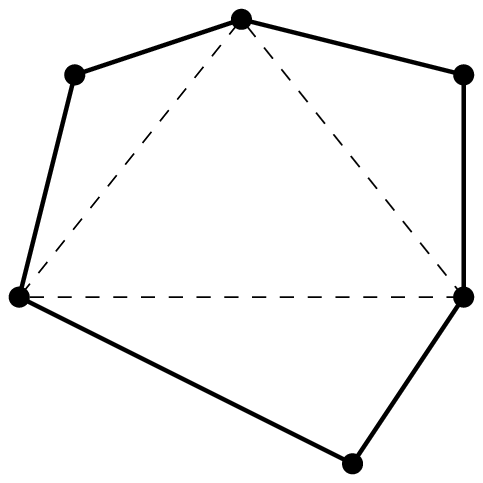}
    \caption{}
    \label{fg:polygonQualityA}
  \end{subfigure}
\hspace{.03\textwidth}
  \begin{subfigure}[b]{0.18\textwidth}
    \includegraphics[height=\textwidth]{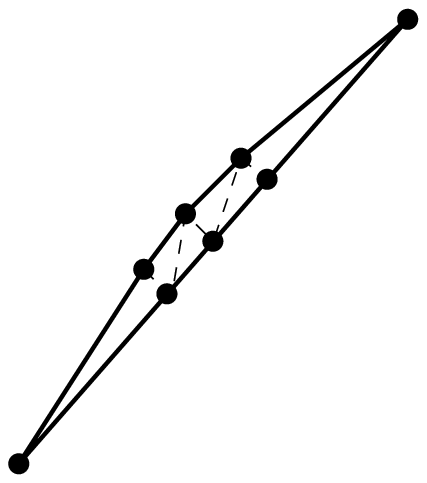}
    \caption{}
    \label{fg:polygonQualityB}
  \end{subfigure}
  \begin{subfigure}[b]{0.18\textwidth}
    \includegraphics[height=\textwidth]{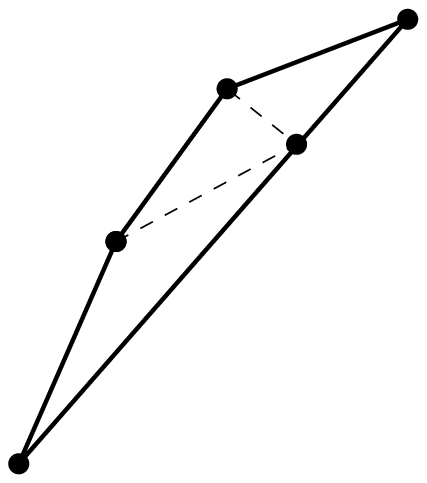}
    \caption{}
    \label{fg:polygonQualityC}
  \end{subfigure}
  \begin{subfigure}[b]{0.18\textwidth}
    \includegraphics[height=\textwidth]{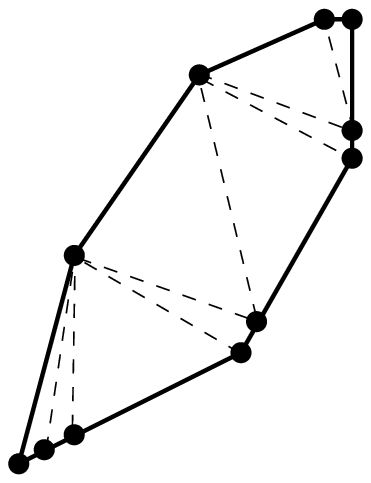}
    \caption{}
    \label{fg:polygonQualityD}
  \end{subfigure}
     \end{center}
\caption{A number of convex polygons with different shape quality characteristics, some of which are distinct from the three classes of triangles:
(A) a polygon with modest aspect ratio and no short edges; 
(B) a very narrow polygon that does not contain any short edges or large angles in the Delaunay triangulation;
(C) a polygon with a large angle in the constrained Delaunay triangulation;
(D) a polygon with good aspect ratio but some very short edges and small angles in the Delaunay triangulation.
The Delaunay triangulation of each polygon is shown. 
(Note: since the polygons are convex, no boundary edges need to be constrained.)
Even if it is part of a larger linear segment (as occurs in (B), (C), and (D)), each segment between vertices shown is treated as an independent edge for the purpose of the upcoming construction (\ref{eq:lineachedge}). 
 }\label{fg:polygonQuality}
\end{figure}

In this paper, we broaden the analysis of harmonic coordinates from \cite{GRB2011} to non-bounded aspect ratio polygons and to three-dimensional polyhedra. 
We connect error in the harmonic coordinates to errors in a piecewise linear approximation for certain triangulations/tetrahedralizations of the polytope. 
Despite the fact that piecewise linear interpolation with respect to a triangulation of the domain and interpolation by harmonic GBCs appear to be at opposite ends of the accuracy spectrum in our previous error analysis, we show here that for convex polygons, the constrained Delaunay triangulation provides an interpolation error that is essentially as good as the harmonic coordinates. 
For non-convex polygons, we establish a similar result connecting large circumradius triangles in the constrained Delaunay triangulation to large interpolation errors, 
although this requires few additional assumptions to avoid some other degenerate configurations that are only possible in the non-convex setting.
We also analyze harmonic coordinates in three dimensions and establish error estimates for bounded aspect ratio polyhedra.
Features of polyhedra that impact interpolation quality are described in both the convex and non-convex cases.
Interpolation on convex polyhedra is analogous to convex polygons although due to sliver tetrahedra, the connection to the Delaunay tetrahedralization cannot be made. 
Rather, poor quality is characterized by vertices being near non-adjacent faces.
For non-convex polyhedra, this characterization does not hold and we give an example demonstrating degenerate limiting geometry with bounded interpolation error.

Section~\ref{sc:tri} begins by stating sharp interpolation results for triangles as a fundamental building block of the analysis.
Section~\ref{sc:gbcs} introduces the key ideas and definitions of generalized barycentric interpolation.
Sections~\ref{sc:twod} and \ref{sc:threed} contain analysis of harmonic coordinates in two and three dimensions, respectively. 
Appendices~\ref{ap:sobolev}, \ref{ap:harmonic}, \ref{ap:bh}, \ref{ap:cdt}, and \ref{ap:tetqual} contain a few well-known results used in our analysis from the theory of Sobolev spaces, partial differential equations, computational geometry, and numerical analysis. 

\section{Barycentric Interpolation on Triangles}\label{sc:tri}

The classical \textit{a priori} analysis of the finite element method involves two steps. 
First the finite element error is bounded (up to some constant factor) by the error in the best possible interpolant in a particular finite element space; this  result is often called Cea's Lemma.
Second, a specific interpolant is constructed and is demonstrated to have a suitably small error in terms of some geometric quantities related to the mesh. 

The most common interpolant selected for linear triangular elements is derived from barycentric coordinates. 
The barycentric functions for a triangle $\displaystyle (\lambda_i)_{i=1}^3$ are the three affine functions taking value $1$ at one vertex of the triangle and taking value $0$ at the other two vertices.
The standard Lagrange interpolant on triangle $T$ is defined in terms of the barycentric coordinates by
\begin{equation}\label{eq:baryinterp}
I_T u = \sum_{i=1}^3 u(\bv_i) \lambda_i
\end{equation}
for any function $u: T\rightarrow \R$ that is smooth enough to admit pointwise values.
In the construction of finite element spaces, this yields two key properties. 
First, the interpolant along a boundary edge of a triangle only depends on the values at the two vertices of the triangle, which allows these functions to be stitched together continuously over the entire mesh.
Second, the space spanned by the barycentric coordinates contains all linear polynomials, which is a crucial piece of the error analysis.

The simplest analysis of the finite element method (widely adopted in the initial exposition of the method in many texts~\cite{Ci02,BS08,EG04,ZT2000}) involves interpolation estimates on triangles under the minimum angle condition:
for minimum angle $\theta>0$, there exists a constant $C_\theta$ such that for all triangles $T$ with minimum angle larger than $\theta$ and all functions $u\in H^2(T)$, 
\begin{equation}\label{eq:minangle}
\hpsn{u-I_T u}{1}{T} \leq C_\theta h_T \hpsn{u}{2}{T}.
\end{equation}
The quantity $h_T$ indicates the length of the longest edge of a triangle $T$.
Under the minimum angle condition, this quantity is equivalent (up to a constant factor) to many other measures of the triangle's ``size,'' such as shortest edge length, inscribed radius, circumradius, diameter, etc. 

Nevertheless, some triangles with arbitrarily small minimum angles exhibit small interpolation errors.
Put differently, only triangles with large angles (near $180^\circ$) yield large interpolation errors. 
Figure~\ref{fg:triangleQuality} depicts the three different classes needed in the refined analysis. 
Under the maximum angle condition, the \textit{a priori} estimate (\ref{eq:minangle}) can be established~\cite{Sy57,BA76,Ja76}, again using the longest edge of the triangle to quantify a triangle's size, which is no longer equivalent to all of the other measures.

Further, a uniform estimate can be established for all triangles that captures the dependence of the constant on the largest angle. 
The most natural form of this estimate avoids any explicit reference to triangle angles by quantifying triangle size using the circumradius.
\begin{theorem}
\label{thm:circumbound}
There exists a constant $C_{\ref{eq:circumbound}}$ such that for all triangles $T$ and all functions $u\in H^2(T)$, 
\begin{equation}\label{eq:circumbound}
\hpsn{u-I_T u}{1}{T} \leq C_{\ref{eq:circumbound}} R_T \hpsn{u}{2}{T},
\end{equation}
where $R_T$ denotes the circumradius of triangle $T$.
\end{theorem}
Throughout the paper, we will use the notation $R_T$ to denote the circumradius of triangle $T$; for an example, see Figure~\ref{fg:delaunayProof} in the appendix.
While some authors have directly connected error estimates over triangles to the circumradius~\cite{Kr92,Ra09,BKK09,KT13,KT13a}, this estimate is also stated in several other essentially equivalent ways in the literature, such as
$h_T/\cos(\alpha_T/2)$ or $h_T/\sin(\alpha_T)$, where $\alpha$ is the largest angle in triangle $T$~\cite{Ja76,GMW99}.

\newpage
\section{Interpolation with Generalized Barycentric Coordinates}
\label{sc:gbcs}

Generalized barycentric coordinates extend some key properties of barycentric coordinates to polygons in a way that allows for analogous finite element spaces to be constructed for polygonal meshes. 
Letting $P$ be a polygon with vertices $(\bv_i)_{i=1}^n$, the functions $(\lambda_i)_{i=1}^n$ are called \textbf{generalized barycentric coordinates} if they satisfy:

\vspace{2mm}
\begin{enumerate}[label=(GBC\arabic*)]
\item\label{gbc:nonneg} {\sc Non-negativity}: $\lambda_i \geq 0$ on $P$;
\item\label{gbc:completeness} {\sc Linear completeness}: For any affine function $L:P\rightarrow \R$, $\displaystyle L = \sum_{i=1}^n L(\bv_i)\lambda_i$;
\item\label{gbc:invariance} {\sc Invariance}: If $F:\R^2 \rightarrow \R^2$ is the composition of rotation, translation and/or uniform scaling operations, then $\lambda_i(\bx) = \lambda_i^F(F(\bx))$, where $\lambda^F$ denotes the barycentric coordinate on $P$;
\item\label{gbc:partition} {\sc Partition of unity}: $\displaystyle \sum_{i=1}^n \lambda_i \equiv 1$;
\item\label{gbc:linearprec} {\sc Linear precision}: $\displaystyle \sum_{i=1}^n \bv_i \lambda_i(\bx) = \bx$;
\item \label{gbc:interp} {\sc Interpolation}: $\lambda_i(\bv_j) = \delta_{ij}$.
\end{enumerate}
\vspace{2mm}

These properties are not independent axioms.
In particular, the last three properties can be derived from the first two. 
We state them together as they identify the complete set of standard properties of generalized barycentric coordinates.

While standard barycentric coordinates on a triangle are unique affine functions that depend only on the location of the vertices, the properties \ref{gbc:nonneg}-\ref{gbc:interp} do not define a unique set of generalized barycentric coordinates in general.
This non-uniqueness has allowed the development of various generalized barycentric coordinate functions tailored to distinct application contexts.
Comparison of the different coordinate types centers around three distinguishing factors:
constraints on the polygonal domain $P$, ease of computation of the coordinates, and the smoothness of the resulting interpolant.
The `interpolant' associated to a set of generalized barycentric coordinates is analogous to the definition in (\ref{eq:baryinterp}), with the sum being taken over all the vertices in $P$.
From this standpoint, we can identify a spectrum of coordinate types ranging from simplicity of implementation to simplicity of analysis.
On one end lie ``triangulation'' coordinates and on the other lie harmonic coordinates, while other types lie somewhere in between.

Triangulation coordinates are the most straightforward generalized barycentric coordinates both conceptually and from the standpoint of computation.
In this case, the polygon is triangulated and a piecewise function is constructed with the barycentric coordinates over the triangles.
This approach comes with some notable downsides.
The triangulation coordinates are not smooth on the interior of the polygon: the interpolant is in $H^1(P)$ but has discontinuous derivatives along the edges of the triangulation.
Moreover, a poor quality triangle can cause large interpolation errors, even if the polygon satisfies certain quality metrics that would make it suitable for other generalized barycentric constructions.
Hence, triangulation coordinates are sensitive to the triangulation selected. 

Harmonic coordinates occupy the opposite end of the spectrum: they are computationally expensive (as they in general have no explicit formula) but are optimal from the perspective of smoothness of the interpolant.
Specifically, the harmonic coordinates are constructed as the solution of the partial differential equation:
\begin{equation}
\label{eq:optpde}
\ds\left\{\begin{array}{rcll}
\Delta \lambda_i& = & 0, & \text{on $P$}, \\
\lambda_i & = & g_i. & \text{on $\p P$}.
\end{array}\right.
\end{equation}
The boundary condition $g_i:\p\Omega\raw\R$ is the piecewise linear function satisfying
\begin{equation}\label{eq:lineachedge}
g_i(\bv_j)=\delta_{ij},\quad g_i \text{ linear on each edge of $\Omega$}.\end{equation}
Since this is a linear PDE, the resulting interpolant can also be characterized as a solution to the differential equation:
\begin{equation}
\label{eq:interpde}
\ds\left\{\begin{array}{rcll}
\Delta \left(I_P u\right)& = & 0 & \text{on $P$}, \\
I_P u & = & g_u & \text{on $\p P$}.
\end{array}\right.
\end{equation}
Here, $g_u$ is the piecewise linear function equal to $u$ at the vertices of $P$. 
We use the same notation for harmonic coordinates as the barycentric interpolant on triangles without confusion, since harmonic coordinates on triangles produce the standard linear interpolant.

These coordinates are optimal in the sense that they minimize the norm of the gradient of the interpolant over all functions satisfying the boundary conditions, that is,
\begin{equation}\label{eq:optimality}
I_P u = \text{argmin} \left\{\hpsn{v}{1}{P} \, : \, v = g_u \,\text{on $\p P$} \right\}.  
\end{equation}
This property, called Dirichlet's Principle (see Appendix~\ref{ap:harmonic}), makes harmonic coordinates ideal for finite element error analysis, since the key estimate in the typical analysis relies on showing that the interpolation operator is bounded from $H^1$ into $H^2$; see~\cite{GRB2011} for more details. 

Many of the other generalized barycentric coordinate types have been developed with essentially the same goal of approximating the harmonic coordinates to avoid large gradients in the interpolant using a construction that does not require the solution of a PDE.
While we only explicitly address harmonic coordinates in this paper, the characterization from (\ref{eq:optimality}) implies that any situation where harmonic coordinates produce large gradients applies to \textit{all} generalized barycentric coordinates. 

\section{Harmonic Coordinates on Polygons}\label{sc:twod}

\subsection{Estimate Based on Triangulation Coordinates}\label{ss:esttri}

We begin with an explicit bound on the harmonic coordinates' interpolation error based on triangulation coordinates. 
This essentially formalizes an expected fact in light of (\ref{eq:optimality}): 
interpolation using harmonic coordinates is no worse than triangulating and using a piecewise linear interpolant on the triangulation.

\begin{theorem}\label{th:triangulationbasedbound}
There exists a constant $C_{\ref{eq:polybound}} > 0$ such that for any polygon $P$, possibly non-convex, all functions $u\in H^2(P)$ and all triangulations $\cT$ of $P$,
\begin{equation}\label{eq:polybound}
\hpsn{u-I_P u}{1}{P} \leq C_{\ref{eq:polybound}} \left(\max_{T\in \cT} R_T \right) \hpsn{u}{2}{P},
\end{equation}
where $R_T$ is the circumradius of a triangle $T$.
\end{theorem}

\begin{proof}
Without loss of generality, we restrict our analysis to a generic diameter one polygon $P$. 
Let $p_u$ denote the linear approximate of $u$ defined by the Bramble-Hilbert lemma; see Theorem~\ref{th:bh} in Appendix~\ref{ap:bh}.
Also noting that linear precision of the harmonic interpolant ensures that $p_u = I_P p_u$, we begin by estimating the error with the triangle inequality:
\begin{align}\label{eq:th41a}
\hpsn{u - I_P u}{1}{P}^2
 & \leq \hpsn{u - p_u}{1}{P}^2 + \hpsn{p_u - I_P u}{1}{P}^2
   = \hpsn{u - p_u}{1}{P}^2 + \hpsn{I_P\left(p_u -  u\right)}{1}{P}^2.
\end{align}

Next we appeal to Dirichlet's principle (see Theorem~\ref{thm:dir-prin} from Appendix~\ref{ap:harmonic}) and observe that the $H^1$-seminorm of the harmonic interpolant is dominated by that of any piecewise linear interpolant that satisfies the same boundary conditions. 
Thus, for any triangulation $\cT$ of $P$, 
\begin{align}\label{eq:th41b}
\hpsn{I_P\left(p_u -  u\right)}{1}{P}^2 \leq \sum_{T\in\cT}\hpsn{I_T\left(p_u -  u\right)}{1}{T}^2.
\end{align}
Combining (\ref{eq:th41a}) and (\ref{eq:th41b}) gives
\begin{align*}
\hpsn{u - I_P u}{1}{P}^2 & \leq \hpsn{u - p_u}{1}{P}^2 + \sum_{T\in\cT}\hpsn{I_T\left(p_u -  u\right)}{1}{T}^2.
\end{align*}

Now while $p_u$ has been selected to be near $u$, it is not directly connected to $I_Tu$. 
So we insert $u$ into the estimate of that term and again apply the triangle inequality:
\begin{align*}
\hpsn{u - I_P u}{1}{P}^2 & \leq \hpsn{u - p_u}{1}{P}^2 + \sum_{T\in\cT}\left[ \hpsn{p_u -  u}{1}{T}^2 + \hpsn{u - I_Tu}{1}{T}^2\right]\\
 & = 2\hpsn{u - p_u}{1}{P}^2 + \sum_{T\in\cT}\hpsn{u - I_Tu}{1}{T}^2.
\end{align*}
The first term is bounded by its construction from the Bramble-Hilbert lemma;
note that the constant $C_{\ref{eq:bh}}$ below comes from Theorem~\ref{th:bh}.
$I_Tu$ is the piecewise linear interpolant of $u$ on triangle $T$ and thus its interpolation error in the second term is bounded by Theorem~\ref{thm:circumbound}.
\begin{align*}
\hpsn{u - I_P u}{1}{P}^2 & = 2 C_{\ref{eq:bh}} \hpsn{u}{2}{P}  + \sum_{T\in\cT}C_{\ref{eq:circumbound}} R_T\hpsn{u}{2}{T}^2\\
 & \leq \left[2 C_{\ref{eq:bh}} + C_{\ref{eq:circumbound}} \max_{T\in\cT} R_T \right] \hpsn{u}{2}{P}.
\end{align*}
Since $\max_{T\in\cT} R_T \geq 1/2$ for any unit diameter polygon, (\ref{eq:polybound}) holds with $C_{\ref{eq:polybound}} := 4 C_{\ref{eq:bh}} + C_{\ref{eq:circumbound}}$.
\end{proof}

\begin{figure}[ht]
     \begin{center}
\includegraphics[width=.3\textwidth]{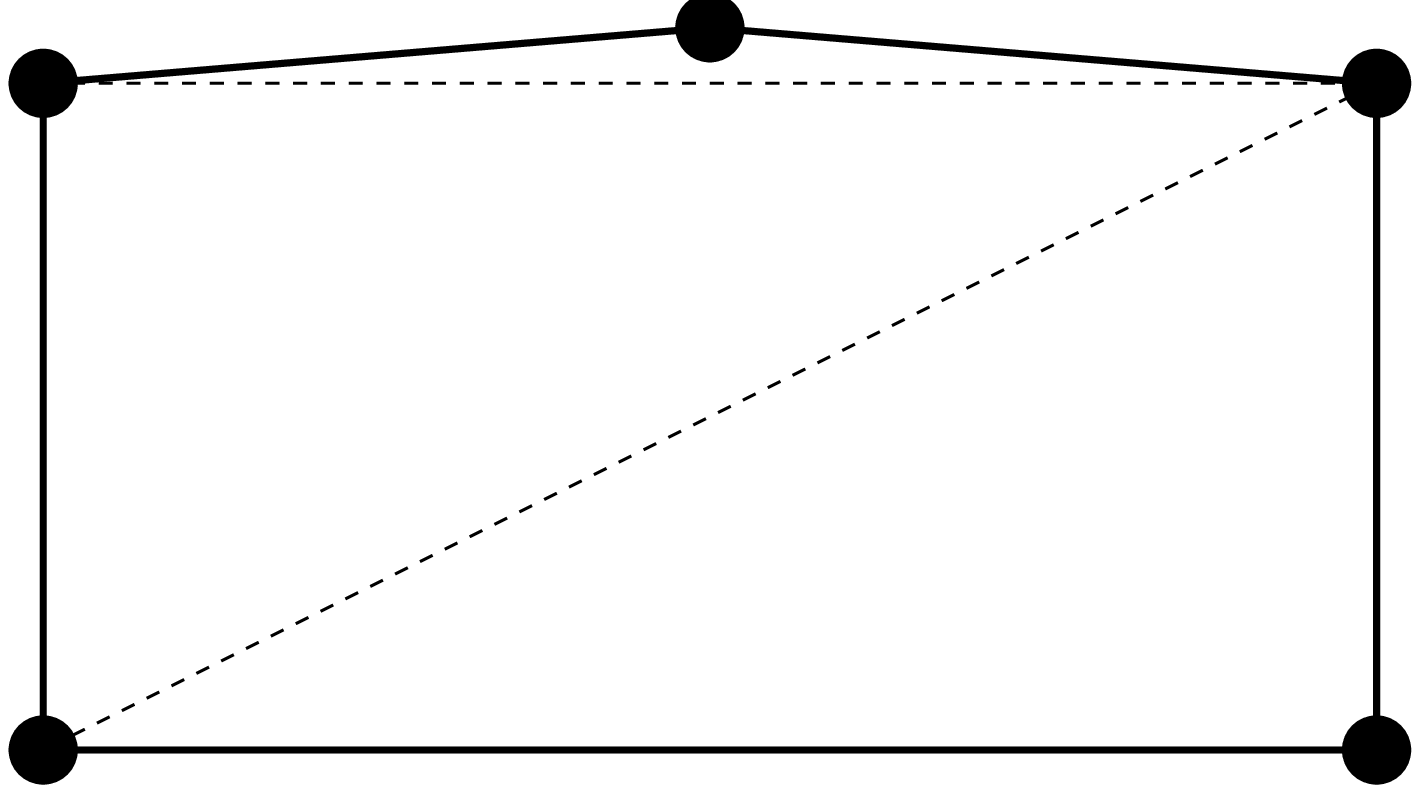}
\hspace{.05\textwidth}
\includegraphics[width=.3\textwidth]{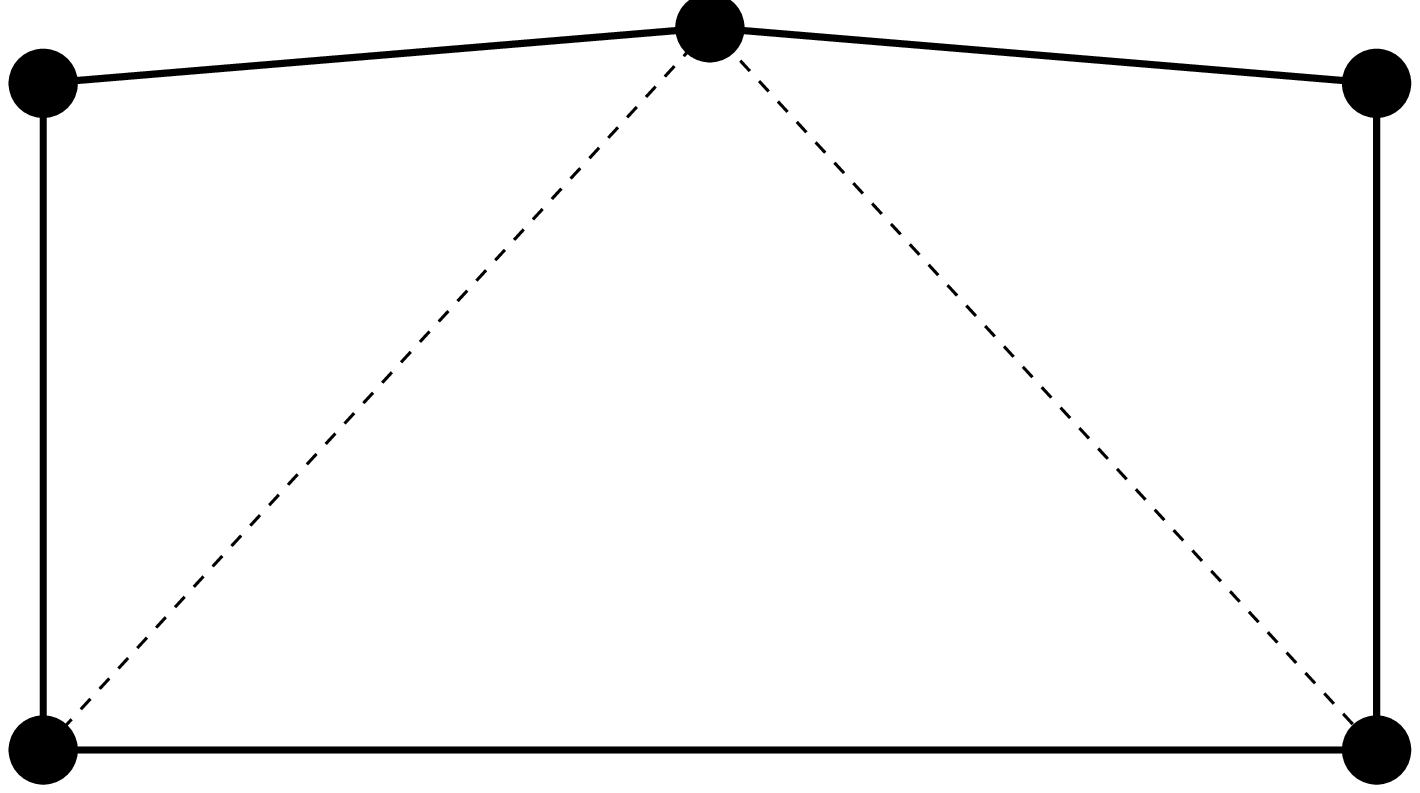}
     \end{center}
\caption{Triangulations of a polygon that leads to a poor (left) and good (right) interpolation error estimate. The constrained Delaunay triangulation on the right is associated with bounded interpolation error.}\label{fg:goodAndBadTriangulation}
\end{figure}

The result just proved begs the question of how to find a triangulation that gives the sharpest estimate.
For some polygons with interior angles near $180^\circ$, certain triangulations may yield arbitrarily poor interpolation properties even though other triangulations yield much smaller errors.
See Figure~\ref{fg:goodAndBadTriangulation} for an example.
In attempt to find the best triangulation for a given polygon, a natural starting point is the Delaunay triangulation, which is based on a geometric criterion (the empty circumball property) that tends to avoid large circumradii triangles.
Since we expect the boundary of the potentially non-convex polygon to be in the triangulation, we are required to use a constrained Delaunay triangulation~\cite{Ch89,Sh08i}, which is a generalization of the Delaunay triangulation that allows for required segments in the triangulation.

Delaunay and constrained Delaunay triangulations have many optimality properties.
Most well known, the Delaunay triangulation maximizes minimum angle~\cite{La77,Si78}, maximizes mean inradius~\cite{La94}, and minimizes $L^p$ error for the function $\vsn{\bx}^2$~\cite{DS89,Ri92,CX04}. 
The most closely connected property is shown in \cite{Ri90, Po92}: the Delaunay triangulation minimizes the $H^1$-semi norm of the interpolant, but not necessarily the error. 
Further, \cite{Ri90a} contains some explicit connections between Delaunay triangulation and the interpolation of harmonic functions.

Finally, we emphasize that Theorem~\ref{th:triangulationbasedbound} applies to all planar polygons, including non-convex and non-bounded aspect ratio polygons that were not included in our original analysis of harmonic coordinates~\cite{GRB2011}. 
This is particularly useful as there are various kinds of polygon degeneracy that do not occur when only considering triangles.
See Figure~\ref{fg:polygonQuality} for some examples.
As a result, fully understanding interpolation estimates for harmonic coordinates on all polygons is reduced to two issues:
(1) showing that the constrained Delaunay triangulation (CDT) for a particular polygon has maximum circumradius bounded by the polygon diameter, or 
(2) demonstrating that the error is in fact large when the CDT contains a large circumradius triangle.

For convex polygons, we will show that this is an essentially complete way of characterizing the interpolation error of harmonic coordinates via the interpolation error of triangulation-based coordinates, excluding some specific families of degenerate polygons.
Non-convex polygons lead to a similar conclusion, although the set of excluded degenerate families is somewhat larger in that setting.

\subsection{Sharpness: Convex Polygons}\label{ss:convex}

In certain situations, harmonic coordinates cannot improve upon interpolation via triangulation coordinates on the constrained Delaunay triangulation.
To make this connection explicit, we begin with a fact connecting large circumradii in the constrained Delaunay triangulation with obtuse triangles involving an edge along the boundary of the polygon. 

\begin{lemma}
Let $P$ be a polygon and let $R_* = \max_{T \in \cT} R_T$, where $\cT$ is  the set of triangles in the constrained Delaunay triangulation of $P$. If $R_\ast > \diam(P)$, then there is a triangle $T_\ast \in \cT$ such that (i) $R_{T_\ast} = R_\ast$ and (ii) the longest edge of $T_\ast$ is on the boundary of $P$. 
\end{lemma}
\begin{proof}
Suppose $R_\ast > \diam(P)$.
Let $T\in \cT$ such that $R_T = R_\ast$. 
If $T$ is acute, then the circumcenter of $T$ lies inside the triangle and thus $\diam(P)\geq R_\ast$, a contradiction.
So $T$ is obtuse.
Now suppose that the longest edge of $T$ is not on the boundary of $P$ and thus is not a constrained edge.
In this case, the longest edge of $T$ is shared by an adjacent triangle $T'$ whose circumradius is at least as large as $R_\ast$ by the Delaunay property; see Figure~\ref{fg:delaunayProof} and~Lemma~\ref{lm:delaunayProp} in Appendix~\ref{ap:cdt}.
Since $T$ has a maximal circumradius by definition, the two circumradii must be the same.
 
We can continue traversing the triangulation by this logic until we reach a triangle along the boundary of $P$. 
This walk will never return to a previously reached triangle, since (1) all triangles encountered must be obtuse and (2) two adjacent obtuse triangles cannot share the same longest edge in a Delaunay triangulation.
Thus the length of the longest edge of the triangles in the walk is strictly increasing. 
Since there are a finite number of triangles, the walk must terminate in a triangle satisfying the conditions specified.
\end{proof}

Next, we demonstrate that the harmonic coordinate interpolant $I_Pu$ will not allow the standard interpolation error estimate over any family of convex polygons with degeneracies that are essentially similar to a family of triangles with arbitrarily large angles.
To this end, we consider a sequence of convex polygons, where one vertex approaches the interior of a non-incident edge and demonstrate that no constant can satisfy the error estimate.
Let $\left(P_i, \bv_i, e_i\right)_{i=1}^{\infty}$ be a sequence of tuples where $P_i$ is a convex polygons, $\bv_i$ is a vertex of $P_i$, and $e_i$ is an edge of $P_i$ not incident to $\bv_i$. 
We assume that this sequence satisfies the following assumptions:
\begin{enumerate}[label=(A\arabic*),series=polygonassumptions]
\item\label{as:diam1} $\diam (P_i) = 1$;
\item\label{as:xaxis} $e_i$ lies on the $x$-axis;
\item\label{as:yaxis} $\bv_i$ lies on the positive $y$-axis;
\item\label{as:slope} $P_i$ lies below a line through $\bv_i$ with non-negative slope;
\item\label{as:nondegen} There exists a $c_v>0$ such that $\dist(\origin, \bv_e) > c_v$, where $\bv_e$ is either endpoint of $e_i$ and $\origin$ is the origin;
\item\label{as:smalldist} $\dist(\bv_i, e_i) \rightarrow 0$ as $i\rightarrow \infty$. 
\end{enumerate}

Assumptions \ref{as:diam1}-\ref{as:slope} can be taken without loss of generality via translation, scaling, rotation and reflection of the domain as needed.
Note that after \ref{as:diam1}-\ref{as:yaxis} have been accommodated, the convexity of $P_i$ ensures that $P_i$ lies below some line through $v_i$.
If that line has negative slope, reflect $P_i$ across the $y$-axis to accommodate assumption \ref{as:slope}.

Assumptions \ref{as:nondegen} and \ref{as:smalldist} are the properties that ensure the polygons are degenerating in a way incompatible with the desired interpolation error estimates.
For instance, the assumptions allow the case of a sequence of triangles with one angle approaching $180^\circ$ (i.e. an unbounded circumradius) by setting $P_i$ to be the triangle with vertices $(-1/2, 0)$, $(1/2, 0)$, and $(0,v_i)$ with $0<v_i<\sqrt 3/2$, and $v_i\raw 0$ as $i\raw\infty$.
The generic configuration specified by \ref{as:diam1}-\ref{as:smalldist} is depicted in Figure~\ref{fg:convexCase}.
For any sequence of polygons satisfying our assumptions, we now show that the interpolation error cannot be bounded by the $H^2$-semi-norm for a particular function.

\begin{lemma}\label{lm:convexcase}
For a sequence of polygon-vertex-edge tuples $\left(P_i, \bv_i, e_i\right)_{i=1}^{\infty}$, satisfying \ref{as:diam1}-\ref{as:smalldist},
\[
\lim_{i\rightarrow \infty} \frac{ \hpsn{u-I_{P_i} u}{1}{P_i} }{ \hpsn{u}{2}{P_i} } = \infty,
\]
for the function $u(x,y) = x^2$.
\end{lemma}

\begin{proof}
First note that $\hpsn{u}{2}{P_i}$ is trivially bounded by a constant:
\begin{equation}\label{eg:denombound}
\hpsn{u}{2}{P_i}^2 = \int_{P_i} 2^2 = 4 \vsn{P_i} \leq \pi, 
\end{equation}
where the final inequality follows from the isodiametric inequality \cite[p.~69]{EG91}.

\begin{figure}[ht]
\[\includegraphics[width=.5\textwidth]{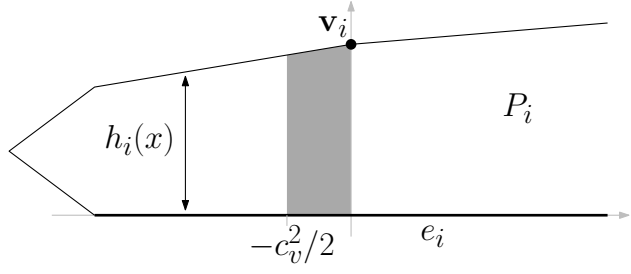}\]
\caption{Generic configuration of a convex polygon of diameter 1 assumed for Lemma~\ref{lm:convexcase}.  For $u(x,y):=x^2$, it is shown that integrating the partial derivative with respect to $y$ of the interpolant of $u$ over the shaded region grows without bound as $d(\bv_i,e_i)\raw 0$.}\label{fg:convexCase}
\end{figure}

It remains to be shown that $\hpsn{u - I_{P_i} u }{1}{P_i}$ grows without bound as $i\rightarrow \infty$. 
We will show that the partial derivative with respect to $y$ blows up near the origin.
Let $h_i(x)$ denote the height of $P_i$ at any value of $x$ along edge $e_i$, as shown in Figure~\ref{fg:convexCase}. 
The convexity of $P_i$, along with assumptions \ref{as:slope} and \ref{as:nondegen}, ensures that  the region bounded by $-c_v\leq x\leq 0$ and $0\leq y\leq h_i(x)$ is a subset of $P_i$.
We will integrate over the even smaller domain, where $-c_v^2/2\leq x\leq 0$; note that assumption \ref{as:diam1} ensures $c_v<1$.
For $u=x^2$, $\frac{\partial}{\partial y} u = 0$ and thus
\begin{align}
\hpsn{u - I_{P_i} u }{1}{P_i}^2 & 
\geq \int_{P_i} \left(\frac{\partial}{\partial y}\left(u - I_{P_i} u \right) \right)^2  
= \int_{P_i} \left(\frac{\partial}{\partial y} I_{P_i} u \right)^2  
\geq \int_{-c_v^2/2}^0 \int_0^{h_i(x)} \left(\frac{\partial}{\partial y}  I_{P_i} u(x,y)\right)^2 \dd y \dd x . \label{eq:subsetintegral}
\end{align}
The inner integral is the $H^1$-semi-norm of a harmonic function, so Dirichlet's principle provides a lower bound.
Using Corollary~\ref{cor:minh1} from Appendix~\ref{ap:harmonic}, we have
\begin{equation}
\label{eq:inn-int-bd}
\int_0^{h_i(x)} \left(\frac{\partial}{\partial y}  I_{P_i} u(x,y)\right)^2 \dd y \geq \frac{\left(I_{P_i} u(x,h_i(x)) - I_{P_i} u(x,0) \right)^2}{h_i(x)}
\end{equation}
We next show that $I_{P_i} u(x,0)-I_{P_i} u(x,h_i(x))\geq c_v^2/2$.
First note that by assumptions \ref{as:yaxis} and \ref{as:nondegen}, edge $e_i$ contains the interval $[-c_v,c_v]$ lying on the $x$-axis.
Since $I_{P_i} u(x,0)$ linearly interpolates the values of $u$ at the endpoints of edge $e_i$ and $u=x^2$ by hypothesis, we have
\begin{equation}\label{eq:interplower}
I_{P_i} u(x,0) \geq c_v^2.
\end{equation} 

\begin{figure}[ht]
\[\includegraphics[width=.25\textwidth]{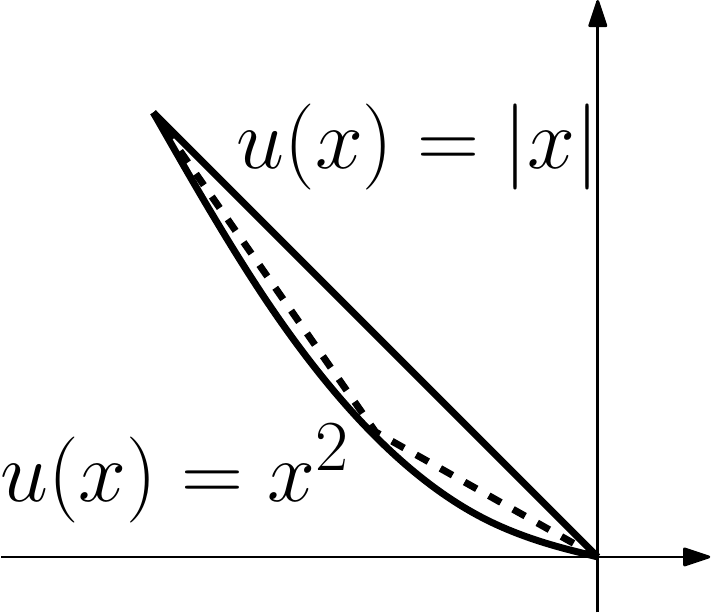}\]
\caption{On $[-1,0]$, any piecewise linear interpolant of $x^2$ lies above $x^2$ and below $|x|$.
This simple fact provides an estimate for $I_{P_i} u(x,h_i(x))$ in the notation of in the proof of Lemma~\ref{lm:convexcase}.}
\label{fg:convexityFigure}
\end{figure}

\noindent
Next, note that for $x\in [-c_v, 0]$, the one-variable function $I_{P_i} u(x,h_i(x))$ reads off the data along a portion of the boundary of $P_i$.
In particular, it is a piecewise linear approximation of $x^2$ that equals $x^2$ at the $x$ values corresponding to the $x$-coordinate values of the vertices of $P_i$.
We have $-1\leq -c_v$ by assumption \ref{as:diam1} and $I_{P_i} u(0,h(0))=0$ since $u=x^2$.  
Thus $x^2 \leq I_{P_i} u(x,h_i(x)) \leq \vsn{x}$ as shown in Figure~\ref{fg:convexityFigure}.
We conclude that for $-c_v^2/2 \leq x \leq 0$, 
\begin{equation}\label{eq:interpupper}
I_{P_i} u(x,h_i(x)) \leq c_v^2/2.
\end{equation}
By (\ref{eq:interplower}) and (\ref{eq:interpupper}), we have $I_{P_i} u(x,0)-I_{P_i} u(x,h_i(x))\geq c_v^2/2$, as claimed.
By assumption \ref{as:slope}, $h_i(x) \leq \dist(\bv_i, e_i)$ for $x\in [-c_v^2/2, 0]$. 
Thus, the right-hand side of (\ref{eq:inn-int-bd}) is bounded below by $(c_v/2)^2/\dist(\bv_i,e_i)$.
Putting this together with (\ref{eq:subsetintegral}), we have that
\begin{align*}
\hpsn{u - I_{P_i} u }{1}{P_i}^2 
& \geq \int_{-c_v^2/2}^0 \frac{\left(c_v^2/2\right)^2}{\dist(\bv_i,e_i)}dx = \frac{c_v^8}{8\,\dist(\bv_i,e_i)}.
\end{align*}
Finally, using (\ref{eg:denombound}) we have that
\[
\frac{ \hpsn{u - I_{P_i} u }{1}{P_i} }{ \hpsn{u}{2}{P_i} } \geq \frac{c_v^8}{8\pi\,\dist(\bv_i,e_i)},
\]
which grows without bound as $i\rightarrow \infty$ by assumption \ref{as:smalldist}.
\end{proof}

\subsection{Example: A Prototypical Family of Non-convex Polygons}\label{ss:nonconvexexample}

Before attempting to expand the result of Lemma~\ref{lm:convexcase} to the case of non-convex polygons, we will consider a simple example containing representative non-convex geometric configuration.
In the example, we begin by attempting to bound the interpolation error despite the fact that most of the assumptions of Section~\ref{ss:convex} hold. 
The analysis connects success of the interpolation estimate to the existence of a certain discontinuous function on the boundary of the domain in a particular Sobolev space. 
For the convex polygons in question, the particular Sobolev space ($H^{1/2}$ on the domain boundary) does not contain discontinuous functions and thus the estimate fails. 
When considering three-dimensional polyhedra, however, the Sobolev embedding does allow the suitable discontinuous functions and thus this example provides a template for the main result in Section~\ref{ss:nonconvex3d}.

\begin{figure}[ht]
\centering
  \begin{subfigure}[b]{0.4\textwidth}
    \includegraphics[width=\textwidth]{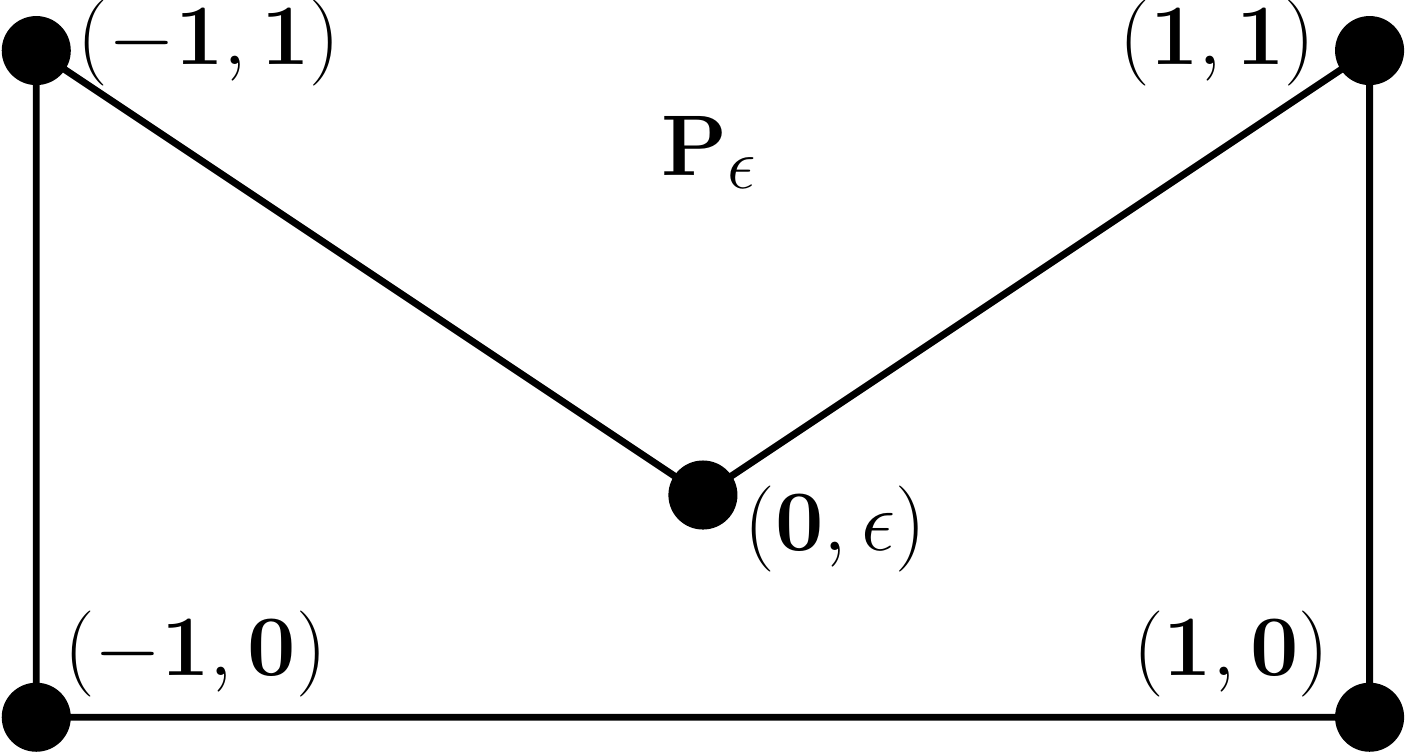}                
    \caption{}
    \label{fg:nonconvexPolygon}
  \end{subfigure}
\hspace{0.08\textwidth}
  \begin{subfigure}[b]{0.25\textwidth}
    \includegraphics[width=\textwidth]{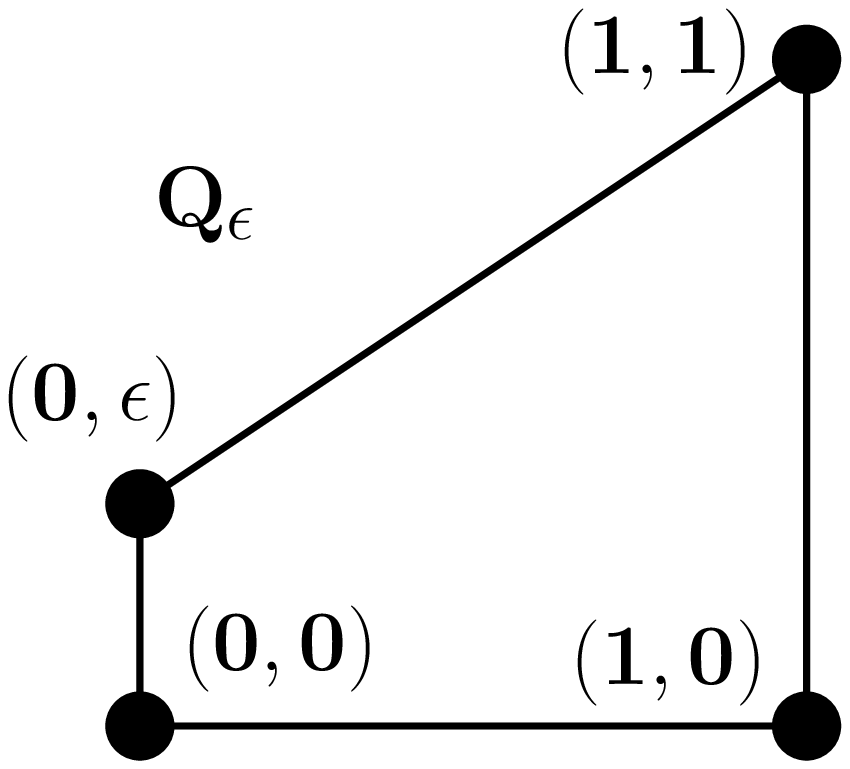}
    \caption{}
    \label{fg:halfDomain}
  \end{subfigure}
\caption{(\textsc{a}) A family of non-convex polygons satisfying assumptions \ref{as:diam1}-\ref{as:yaxis} and \ref{as:nondegen}-\ref{as:smalldist}, which is analyzed in Section~\ref{ss:nonconvexexample}.  (\textsc{b}) The discussion makes use of the domain $Q_\eps$, which is a subset of $P_\eps$ as shown.}
\label{fig:P-and-Q}
\end{figure}

Let $P_\epsilon$ be the five-sided non-convex polygon with vertices at $(-1,0)$, $(-1,1)$, $(1,0)$, $(1,1)$, and $(0,\epsilon)$, with $0<\eps<1$, as depicted in Figure~\ref{fg:nonconvexPolygon}, and consider interpolating the function $u(x,y) = x^2$ on $P_\epsilon$. 
We will use the notation $P_0$ to refer to the limiting domain of the union of two triangles.
Observe that the sequence associated to $\{P_{1/i}\}$ for $i=2,3,4,\ldots$, when scaled by a factor of $1/\sqrt 5$, satisfies assumptions \ref{as:diam1}-\ref{as:yaxis} and \ref{as:nondegen}-\ref{as:smalldist}.
The scaling by the fixed value $1/\sqrt 5$ does not affect the subsequent argument, so we focus on estimating $\hpsn{u-I_{P_\eps} u}{1}{P_\eps}/\hpsn{u}{2}{P_\eps}$ as $\eps\raw 0$ without explicitly scaling the domain to diameter one. 
Note also that the constrained Delaunay triangulation of $P_\eps$ contains a triangle $T_\epsilon$ with vertices $(-1,0)$, $(1,0)$ and $(0,\epsilon)$ for which the circumradius grows as $\epsilon \rightarrow 0$, so Theorem~\ref{th:triangulationbasedbound} cannot be employed.

Observe that $\hpsn{u}{2}{P_\epsilon}$ is bounded below uniformly and $\hpsn{u}{1}{P_\epsilon}$ is bounded above uniformly:
\begin{align}
\hpsn{u}{2}{P_\epsilon} ^2  \geq \hpsn{u}{2}{P_0}^2 = \int_{P_0} 2 ^2 = 4;\label{eq:h2lower}\\
\hpsn{u}{1}{P_\epsilon} ^2  \leq \hpsn{u}{1
}{P_1}^2 = \int_{P_1} \left(2x \right)^2 = 8/3.\label{eq:h1upper}
\end{align}
By the triangle inequality 
\begin{equation}\label{eq:triineq}
\hpsn{u - I_{P_\epsilon} u}{1}{P_\epsilon} \leq \hpsn{u}{1}{P_\epsilon} + \hpsn{ I_{P_\epsilon} u}{1}{P_\epsilon},
\end{equation}
thus the remaining term to estimate is $\hpsn{I_{P_\epsilon} u}{1}{P_\epsilon}$.
Standard approaches toward this goal (e.g.\ the Sobolev embedding theorem, Morrey's inequality) cannot be established uniformly in a trivial way due to the non-locality of the boundary, where the domain becomes arbitrarily narrow near the origin.
To skirt this difficulty, we work with the subdomain $Q_\epsilon := \{ (x,y) \in P_\epsilon \,|\, x \geq 0\}$, as shown in Figure~\ref{fg:halfDomain}.

Since $u(x,y)=x^2$, the boundary condition $g_\eps:\p P_\eps\raw\R$ is identically 1 along the horizontal and vertical edges of $P_\eps$ and decays linearly from 1 at $(\pm 1, 1)$ to zero at $(0,\eps)$ along the two diagonal edges.
By Dirichlet's principle (see Theorem~\ref{thm:dir-prin} from Appendix~\ref{ap:harmonic}), $I_{P_\epsilon} u$ has minimal $H^1$-semi-norm among all functions satisfying this boundary condition.
By the symmetry of the domain and the boundary condition, Dirichlet's principle also holds over the subdomain $Q_\epsilon$, i.e.,
\begin{equation}
\label{eq:optimality2}
\hpsn{I_{P_\epsilon} u}{1}{P_\epsilon}^2 = 2\hpsn{I_{P_\epsilon} u}{1}{Q_\epsilon}^2  \leq 2 \hpsn{v}{1}{Q_\epsilon}^2,
\end{equation}
for any $v\in H^1(Q_\epsilon)$ that satisfies the same boundary condition as $I_{P_\epsilon} u$ on all the boundary segments shared between $P_\epsilon$ and $Q_\epsilon$.
In particular, (\ref{eq:optimality2}) applies in the case $v:={\rm Tr}^{-1}g_\epsilon$, where ${\rm Tr}^{-1}$ denotes the right-continuous inverse of the trace operator for the domain $Q_\eps$. 
Note that the trace operator and its inverse depend on the domain, i.e.\ on $\eps$, but we omit this explicit dependence from our notation as dependence on $\eps$ is already indicated by the boundary data $g_\eps$.
We can now estimate
\begin{equation}\label{eq:invtracebound}
\hpsn{I_{P_\epsilon} u}{1}{P_\epsilon}^2 \leq 2 \hpsn{{\rm Tr}^{-1}g_\epsilon}{1}{Q_\epsilon}^2 \leq C_{\ref{eq:invtracebound}}\hpn{g_\epsilon}{1/2}{\partial Q_\epsilon}^2.
\end{equation}
For Lipschitz domains, the trace operator has a well-defined continuous, bounded inverse, meaning that for any fixed $\epsilon$, the constant $C$ above exists; see Theorem~\ref{th:trace} in Appendix~\ref{ap:sobolev} for more details.

Before analyzing $\hpn{g_\epsilon}{1/2}{\partial Q_\epsilon}$, we must argue that the constant $C_{\ref{eq:invtracebound}}$ is independent of $\epsilon$. 
The typical technique in the analysis of Sobolev-type properties on Lipschitz domains, e.g.\ extension or trace theorems, involves covering the domain by a finite number of overlapping patches, constructing a partition of unity supported by the patches, and then on each patch intersecting and ``flattening'' the boundary via a Lipschitz map.
The bound on the inverse trace operator depends on the domain only by the parameters of this decomposition: the number of mutually overlapping patches, the behavior of the partition of unity, and the Lipschitz constant of the maps for flattening the boundary. 
For $Q_\eps$, this construction can be performed with a set of patches and partition of unity that is independent of $\eps$.
Further, the Lipschitz constant of the flattening maps can also be bounded independent of $\eps$ as shown in Lemma~\ref{lm:flattening} in Appendix~\ref{ap:sobolev}.
Thus the constant $C_{\ref{eq:invtracebound}}$ can be taken independent of $\eps$.

\begin{figure}[ht]
     \begin{center}
\includegraphics[width=.3\textwidth]{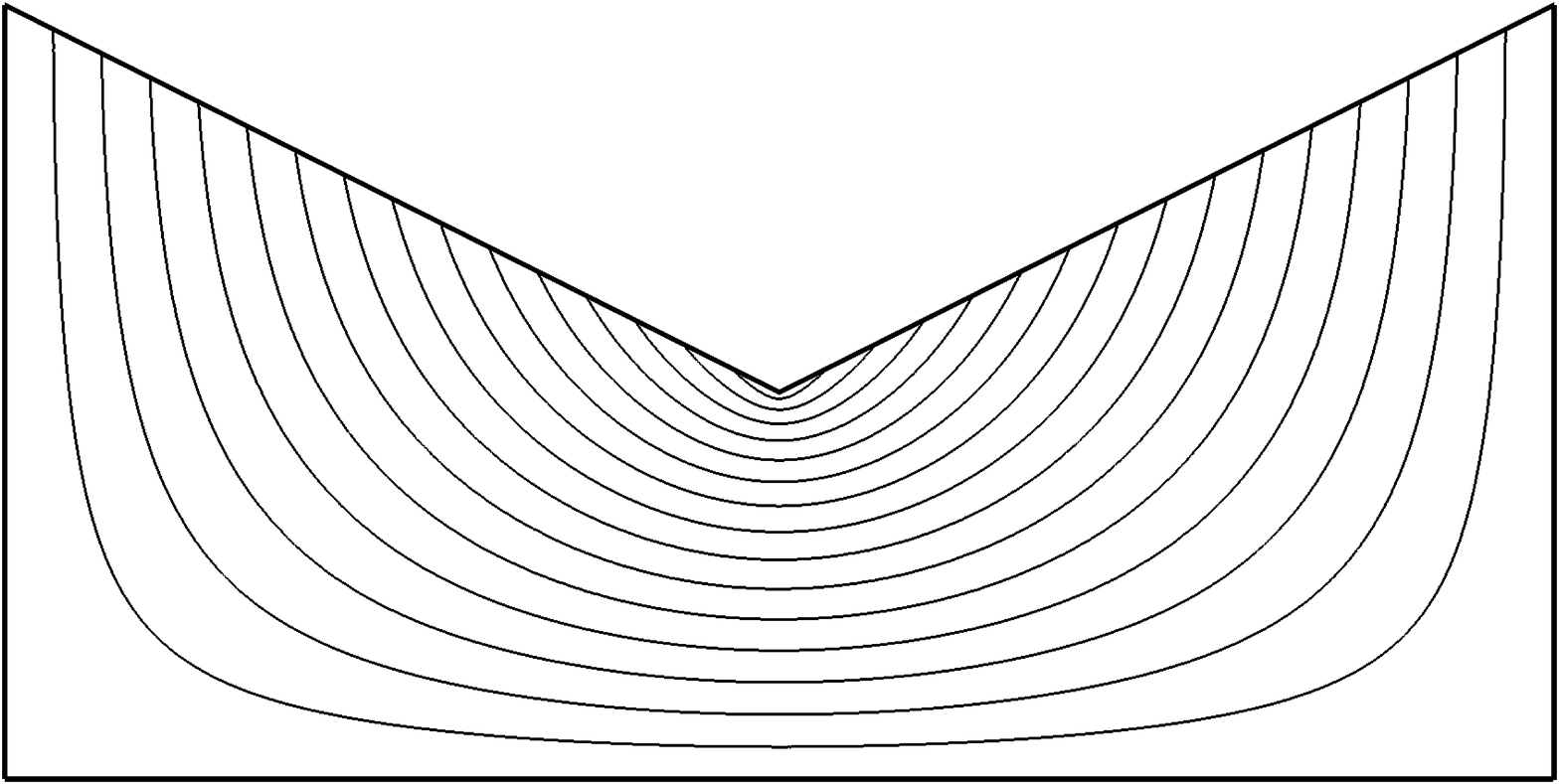}
\includegraphics[width=.3\textwidth]{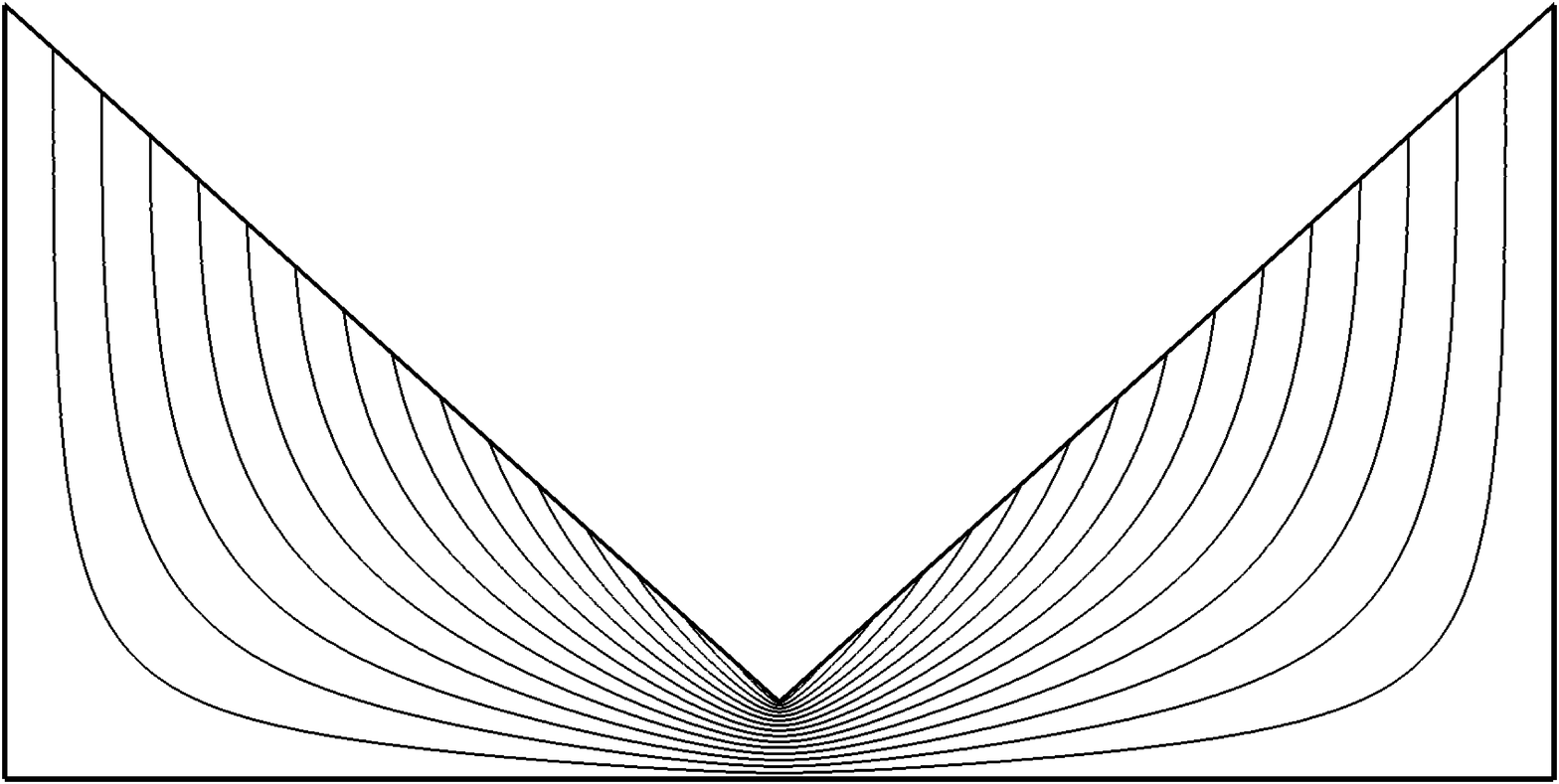}
\includegraphics[width=.3\textwidth]{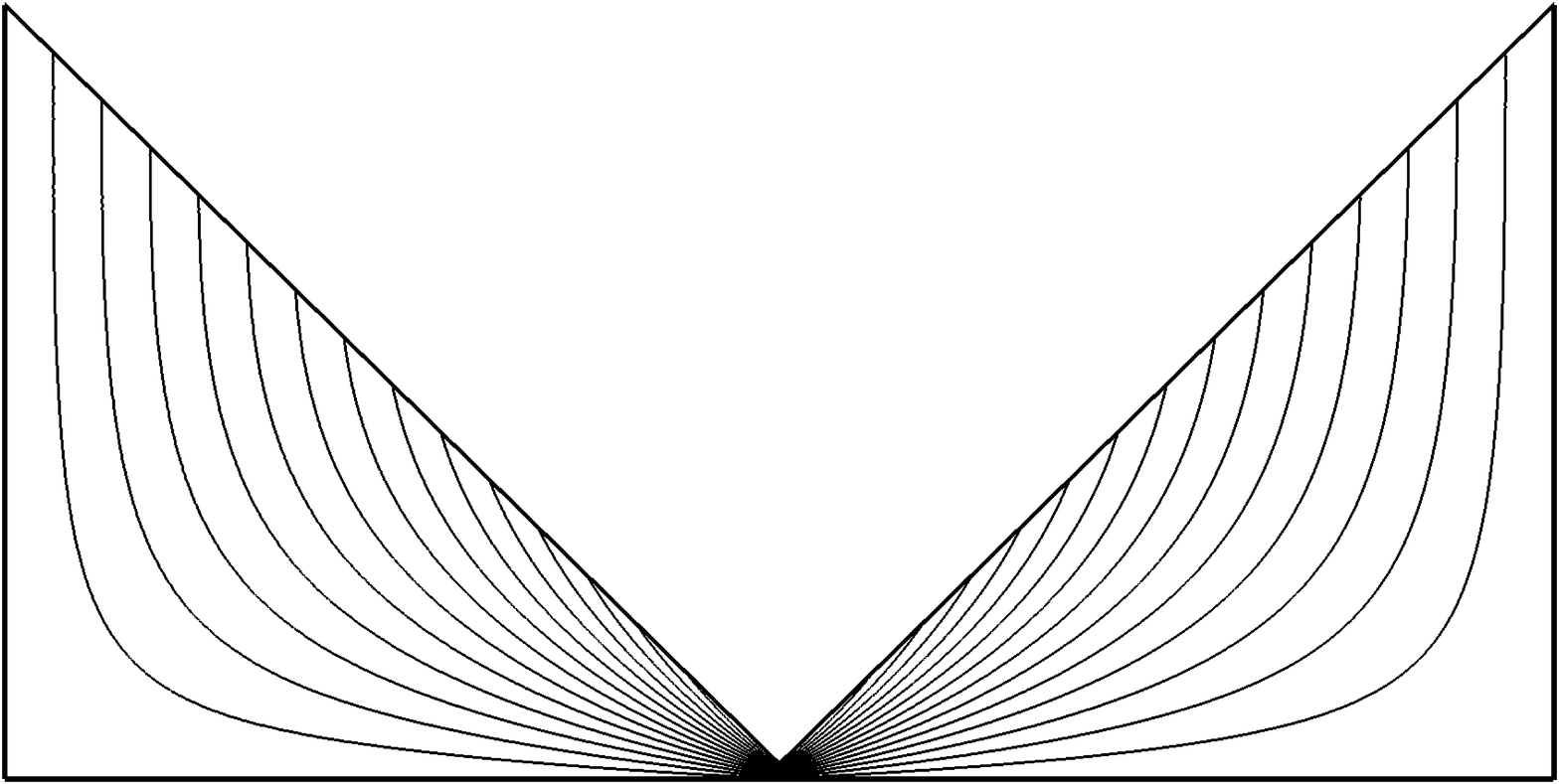}
     \end{center}
\caption{Contour plots of $I_{P_\epsilon}u$ for (left to right) $\epsilon = 0.5$, $\epsilon = 0.1$, and $\epsilon = 0.02$. Large gradients are concentrated near a single point, but, in two dimensions, this is not not enough to cause the interpolation error to be unbounded in the limit.}
\label{fg:contourPlots}
\end{figure}

We examine the limiting behavior of the estimate by considering the limit of $g_0$, i.e.\ the limit of the boundary condition.
Here, $g_0$ is discontinuous on the one-dimensional boundary of the limiting two-triangle domain $Q_0$.
In one dimension, the space $H^{1/2}$ contains only continuous function, as guaranteed by the 
Sobolev embedding theorem or, more specifically, by Morrey's inequality (stated in Theorem~\ref{th:morrey} in Appendix~\ref{ap:sobolev}).
Thus $g_\epsilon$ does not converge in $H^{1/2}(\partial Q_\epsilon)$ as $\eps\raw 0$ and hence $\hpn{g_\epsilon}{1/2}{\partial Q_\epsilon}$ cannot be bounded uniformly.

\begin{remark}\label{rm:sobo}
\em
The Sobolev embedding theorem ensures the continuity of functions in $W^{m,p}$ when $mp \geq d$ where $d$ is the dimension of the space. 
The setting relevant to the example above ($d=1$, $m=1/2$, $p=2$) is the equality case of this inequality.
The case of polyhedra in three dimensions involves combinations that do not satisfy the inequality ($d=2$, $m=1/2$, $p=2$) and the resulting conclusions will be different.
These cases are discussed in Section~\ref{ss:nonconvex3d}. 
Similarly, analyzing harmonic coordinates in $W^{m,p}$ spaces for $p > 2$ can also lead to error bounds, even for polygons; in those cases, however, there are other hurdles, such as the inability to directly appeal to Dirichlet's principle.
\end{remark}

\subsection{Sharpness: Non-convex Polygons}\label{ss:nonconvexresult}
To analyze the sharpness of estimate (\ref{eq:polybound}) for families of non-convex polygons, we begin with the same assumptions that preceded Lemma~\ref{lm:convexcase}, except for \ref{as:slope}, which is essentially tied to the local convexity of the polygon at the vertex in question $\bv_i$. 
This will be replaced by some other restrictions to ensure that non-adjacent entities of the polygon do not interfere with the interaction of vertex $\bv_i$ and non-incident edge $e_i$.

Specifically, we restrict ourselves to families of polygons satisfying \ref{as:diam1}-\ref{as:yaxis}, \ref{as:nondegen}-\ref{as:smalldist}, and the following:
\begin{enumerate}[resume*=polygonassumptions]
\item\label{as:nondegen1} $\dist(\bv_i, \bv_a)> c_v$ for all vertices of $P_i$ where $\bv_i \neq \bv_a$;
\item\label{as:nondegen2} $\dist(\bv_i, e_a) > c_v$ for all edges of $P_i$  other than $e_i$ and the two edges incident to $\bv_i$;
\item\label{as:rightside} the segment between $\bv_i$ and the origin is contained in $P_i$.
\end{enumerate}

\begin{figure}[ht]
     \begin{center}
  \begin{subfigure}[b]{0.18\textwidth}
    \includegraphics[height=\textwidth]{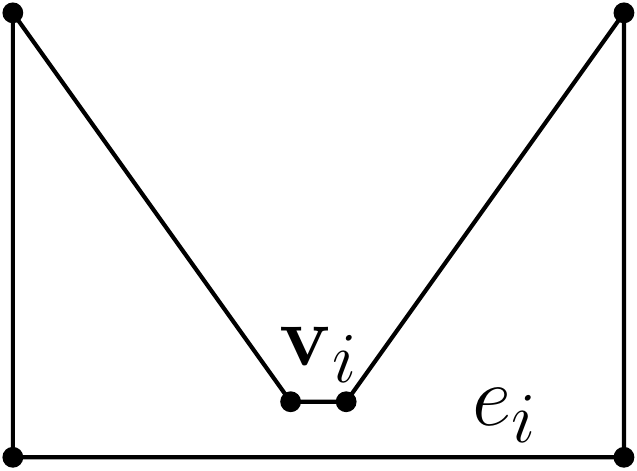}
    \caption{}
    \label{fg:nonconvexRestrictionsA}
  \end{subfigure}
\hspace{.08\textwidth}
  \begin{subfigure}[b]{0.18\textwidth}
    \includegraphics[height=\textwidth]{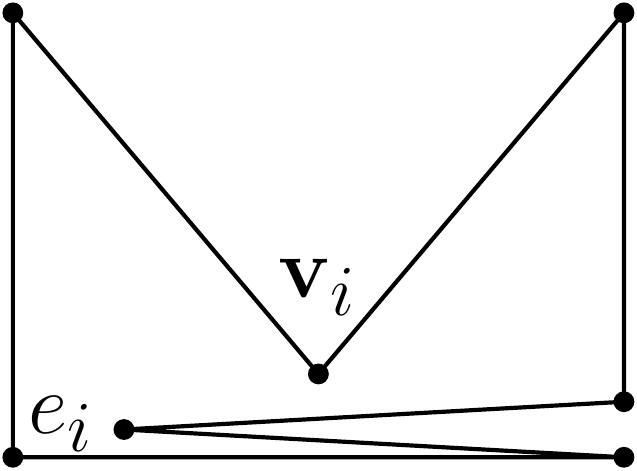}
    \caption{}
    \label{fg:nonconvexRestrictionsB}
  \end{subfigure}
\hspace{.08\textwidth}
  \begin{subfigure}[b]{0.18\textwidth}
    \includegraphics[height=\textwidth]{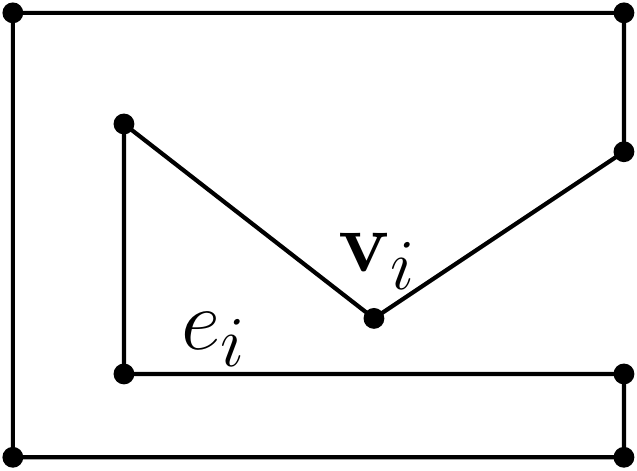}
    \caption{}
    \label{fg:nonconvexRestrictionsC}
  \end{subfigure}
       \end{center}
\caption{Several non-convex polygons that violate various restrictions in Lemma~\ref{lm:nonconvexcasenew}:
(A) a polygon with two nearby vertices, violating \ref{as:nondegen1}; 
(B) a polygon with a vertex near multiple non-incident edges, violating \ref{as:nondegen2}; 
(C) a polygon with a vertex approaching a non-incident edge from the outside of the polygon,violating \ref{as:rightside}. 
 }\label{fg:nonconvexRestrictions}
\end{figure}

In the above assumptions, $c_v$ is the same constant as in \ref{as:nondegen}, although this holds without loss of generality since it can be assumed to be the minimum constant among \ref{as:nondegen}, \ref{as:nondegen1}, and \ref{as:nondegen2}.
For polygons with at least four sides, \ref{as:nondegen2} is strictly stronger than \ref{as:nondegen1} since any vertex near $\bv_i$ also will be connected to an edge that is not incident to $\bv_i$. Still, we include \ref{as:nondegen1} for clarity.
The assumption \ref{as:rightside} serves to exclude situations in which $\bv_i$ approaches $e_i$ from the outside of the $P_i$ which generally does not cause the estimate to blow up; see Figure~\ref{fg:nonconvexRestrictions}.

\begin{lemma}\label{lm:nonconvexcasenew}
For a sequence of polygon-vertex-edge tuples $\left(P_i, \bv_i, e_i\right)_{i=1}^{\infty}$, satisfying \ref{as:diam1}-\ref{as:yaxis} and \ref{as:nondegen}-\ref{as:rightside},
\[
\lim_{i\rightarrow \infty} \frac{ \hpsn{u-I_{P_i} u}{1}{P_i} }{ \hpsn{u}{2}{P_i} } = \infty,
\]
for the function $u(x,y) = x^2$.
\end{lemma}
\begin{proof}
First, we identify some some implications of the assumptions on the domain and, without loss of generality, define a generic case for analysis, which is depicted in Figure~\ref{fg:nonconvexCaseNew}.
Specifically, let $d_{i,1}$ and $d_{i,2}$ be the two edges of $P_i$ incident to vertex $v_i$. 
Without loss of generality, we assume that $d_{i,1}$ lies in the second quadrant, i.e.\ $x < 0$, since reflecting the domain across the $y$-axis does not change the subsequent analysis. 
Moreover, if $d_{i,2}$ is also in the second quadrant, we assume that $d_{i,1}$ is ``below'' $d_{i,2}$, or, more formally, a counter-clockwise sector from $d_{i,1}$ to $d_{i,2}$ is interior to the polygon.
Finally, we restrict our analysis to the case where the line through $d_{i,1}$ has negative slope: if the line has positive slope, then we have the configuration as in the proof of Lemma~\ref{lm:convexcase} and can apply the same argument. 

\begin{figure}[ht]
     \begin{center}
\includegraphics[width=.5\textwidth]{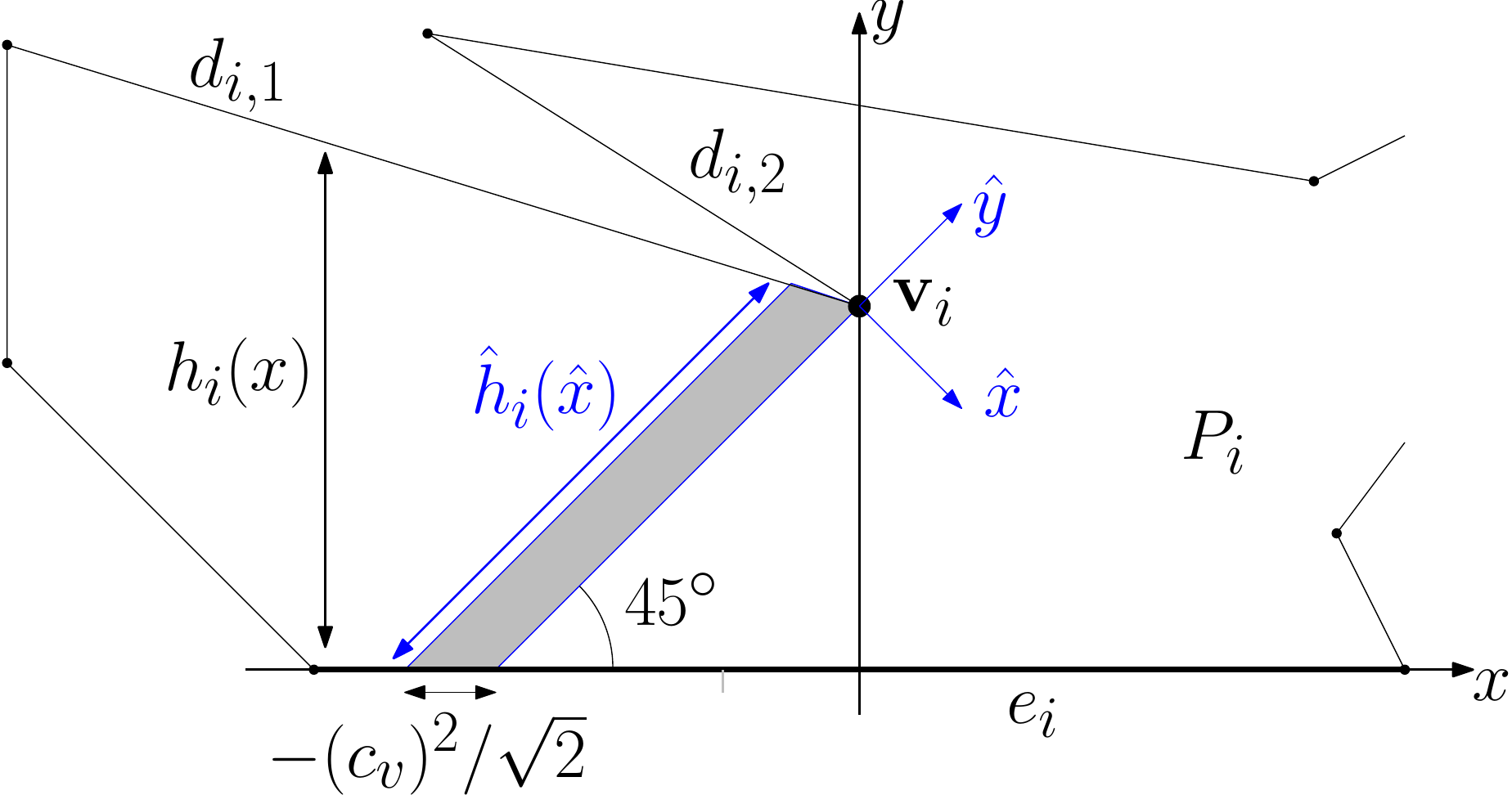}
     \end{center}
\caption{Generic configuration of a polygon assumed for Lemma~\ref{lm:nonconvexcasenew}, showing the standard and rotated coordinate systems $(x,y)$ and $(\hat x, \hat y)$, respectively.  For $u(x,y):=x^2$, it is shown that integrating the partial derivative with respect to $\hat{y}$ of the interpolant of $u$ over the shaded region grows without bound as $d(\bv_i,e_i)\raw 0$.}\label{fg:nonconvexCaseNew}
\end{figure}

As in Lemma~\ref{lm:convexcase}, we will identify a suitable integration subdomain and show that the interpolation estimate must blow up simply looking at the subdomain.
Estimate (\ref{eg:denombound}) still holds, providing a constant upper bound on the denominator: $\hpsn{u}{2}{P_i} \leq \pi$. 
For the numerator, first note that
\begin{align*}
\hpsn{u - I_{P_i} u }{1}{P_i}
\geq \hpsn{I_{P_i} u }{1}{P_i} - \hpsn{u}{1}{P_i}.
\end{align*}
The term $\hpsn{u}{1}{P_i}$ can also be bounded above by a constant, similar to $\hpsn{u}{2}{P_i}$, and thus we only need to show that $\hpsn{I_{P_i} u }{1}{P_i}$ grows without bound as $i\rightarrow \infty$. 

Establishing a lower bound on $\hpsn{I_{P_i} u }{1}{P_i}$ is similar to (\ref{eq:subsetintegral}), but involves analyzing a different integration domain.
To describe this domain, consider an alternative coordinate system, with axes denoted $\hat x$ and $\hat y$, centered at $\bv_i$ and rotated $45^\circ$ clockwise: in other words, the $\hat x$ axis lies on a line of slope $-1$ and the $\hat y$ axis lies on a line of slope $+1$, as shown in Figure~\ref{fg:nonconvexCaseNew}. 
We restrict attention to the region between edges $e_i$ and $d_{i,1}$ in the range $\hat x \in (-c_v/2,0)$, which is the shaded region in Figure~\ref{fg:nonconvexCaseNew}.
Denote this region by $D_i$.
For sufficiently large $i$, assumption \ref{as:smalldist} in conjunction with \ref{as:nondegen1}-\ref{as:rightside} ensure that $D_i \subset P_i$; from here on, we always take $i$ large enough that this containment holds. 
Further, the $(\hat x,\hat y)$ coordinate system has been constructed so that in our generic configuration, with $d_{i,1}$ downward sloping in the second quadrant, the distance between $d_{i,1}$ and $e_i$ is small in the $\hat y$ direction. 
Estimating $\hpsn{I_{P_i} u }{1}{P_i}$ over the domain $D_i$ gives:
\begin{align}
\hpsn{I_{P_i} u }{1}{P_i}^2 & 
\geq \int_{D_i} \left(\frac{\partial}{\partial \hat{y}}\left(I_{P_i} u \right) \right)^2  
\geq \int_{-c_v^2/2}^0 \int_{L_i(\hat x)} \left(\frac{\partial}{\partial \hat y}  I_{P_i} u(\hat x,\hat y)\right)^2 \dd \hat y \dd \hat x, \label{eq:subsetintegralnew}
\end{align}
where $D_i$ represents the integration subdomain associated with polygon $P_i$ and $L_i(\hat x)$ denotes the appropriate slice of the integration domain for a given $\hat x$ value.
To estimate the inner integral, we observe that Dirichlet's principle can be applied as was done in (\ref{eq:inn-int-bd}) and thus, 
\begin{align*}
\int_{L_i(\hat x)} \left(\frac{\partial}{\partial \hat y}  I_{P_i} u(\hat x,\hat y)\right)^2 \dd \hat y 
\geq \frac{\left(I_{P_i} u(\hat x,d_{i,1}(\hat x)) - I_{P_i} u(\hat x,e_i(\hat x)) \right)^2}{\vsn{L_i(\hat x)}},
\end{align*}
where $(\hat x,d_{i,1}(\hat x))$ and $(\hat x,e_{i}(\hat x))$ denote the end points of $L_i(\hat x)$ on $d_{i,1}$ and $e_i(\hat x)$, respectively, and $\vsn{L_i(\hat x)}$ is the length of segment $L_i(\hat x)$. 

As in Lemma~\ref{lm:convexcase}, $u$ along $e_i$ is at least $c_v^2$, and $u$ along $d_{i,1}$ is at most $c_v^2/2$. Thus (\ref{eq:subsetintegralnew}) becomes,
\begin{align}
\hpsn{I_{P_i} u }{1}{P_i}^2 & 
\geq \int_{-c_v^2/2}^0 \frac{c_v^4}{4\vsn{L_i(\hat x)}} \dd \hat x. \label{eq:subsetintegralnew2}
\end{align}
The integration set $L_i(\hat x)$ is contained in the second quadrant and a line of slope $+1$, and thus we can estimate 
\begin{align*}
\vsn{L_i(\hat x)} \leq 2\left(\dist(\bv_i, e_i)/\sqrt{2} - \hat x \right) ,
\end{align*}
noting that $\hat x < 0$ in our area of interest. 
Thus (\ref{eq:subsetintegralnew2}) becomes
\begin{align*}
\hpsn{I_{P_i} u }{1}{P_i}^2 & 
\geq \int_{-c_v^2/2}^0 \frac{c_v^4}{8\left(\dist(\bv_i, e_i)/\sqrt{2} - \hat x \right)} \dd \hat x 
 = \frac{c_v^4}{8}\left(\ln \left(\dist(\bv_i, e_i)/\sqrt{2} + c_v^2/2\right) - 
  \ln \left(\dist(\bv_i, e_i)/\sqrt{2}\right)\right).
\end{align*}
Taking the limit as $i\raw\infty$, the first term on the right approaches a constant while  \ref{as:smalldist} ensures that the second term grows without bound, thus completing the result.
\end{proof}

\section{Harmonic Coordinates on Polyhedra}
\label{sc:threed}

As is typical in computational geometry, going from two to three dimensions introduces additional subtleties and challenges.
On tetrahedra in $\R^3$, as on any $n$-simplex in $\R^n$, barycentric functions are uniquely defined and linear, yielding a similar analysis of their interpolation properties.
On a general polyhedron $P$, however, just defining generalized barycentric coordinates on $\p P$ is not unique or even trivial, since for non-triangular boundary facets, a two-dimensional generalized barycentric coordinate must be employed. 
Once values on $\p P$ are fixed by some means, they are used as boundary conditions to define the harmonic coordinates by the same PDE (\ref{eq:optpde}) from the two-dimensional construction.

In Section~\ref{ss:convexbar}, \emph{we restrict our analysis to polyhedra with triangular facets}, where the GBCs on the boundary are linear functions.
This is the most common setting for existing work on 3D GBCs~\cite{FKR2005,JSW05,JMRGS07,JLW07} although some constructions have been developed and analyzed more generally~\cite{W1996,WSHD2007,HS2008,FGS2013}.
While the result should still hold for a broader class of polyhedra with non-triangular faces, the restricted setting helps to avoid complexity associated with the two-dimensional GBCs on the boundary facets, for which shape-quality and error estimates are already significantly more intricate than in the triangular setting.

In Section~\ref{ss:convex3d}, we consider arbitrary convex polyhedra and demonstrate that the error estimate will fail in any family of polyhedra where vertices are allowed to be arbitrarily close to the interior of non-adjacent facets. 
Finally, in Section~\ref{ss:nonconvex3d}, we show that non-convex polyhedra do not have this property; vertices can approach opposite faces in certain ways such that the error does not blow up. 
Comparisons to the 2D analogues of each result are made in each section.

\subsection{Convex, Bounded Aspect Ratio Polyhedra}\label{ss:convexbar}

First we define the set of polyhedra for which error estimates analogous to the result in \cite{GRB2011} can be established. 
For a convex polyhedron $P$, 
an \textit{insphere} of $P$ is a sphere inscribed in $P$ of maximum radius $\rho(P)$.
Given $P$, the \textit{inradius} $\rho(P)$ is uniquely defined and we fix some particular insphere center $\bc\in P$ arbitrarily if there is not a unique insphere.
The inradius of a two-dimensional facet $F$ of $P$ is defined as in the polygonal case and denoted $\rho(F)$.
The \textit{aspect ratio}, also called the \textit{chunkiness parameter}, is denoted by $\gamma$ and is defined as the ratio of the diameter to the inradius of $P$, i.e.
\[\gamma(P) := \frac{\diam(P)}{\rho(P)}.\]

Given a bound $\gamma_* > 0$ on the aspect ratio of both the polyhedron and each of its faces, let $\cP$ be the set of polyhedra $P$ with the following properties:
\begin{itemize}
\item $P$ is convex;
\item all facets of $P$ are triangles;
\item $\gamma(P) < \gamma_*$; 
\item $\rho(F) > \diam(P)/\gamma_*$ for each facet $F$ of $P$.
\end{itemize}

These restrictions immediately impose a bound on the number of faces and vertices of the polygon as stated in the following lemma. 

\begin{lemma}
There exists $n_* > 0$ depending only upon $\gamma_*$ such that all polygons of $\cP$ contain fewer than $n_*$ triangles and vertices. 
\end{lemma}
\begin{proof}
Without loss of generality (up to affine scaling), we assume that $P$ has diameter $1$. 
Since the filled sphere has maximum surface area over all convex sets of equal diameter~\cite{Bl1915}, the surface area of any $P\in\cP$ is at most $\pi$. 
The inradius of each face of $P$ is at least $\frac{1}{\gamma_*}$, 
so estimating very conservatively, there are at most $\pi \gamma_*^2$ triangles forming the boundary of $P$. 
Euler's formula for the boundary triangulation of the polygon gives the following relationship between the number of vertices $n_v$ and the number of triangles $n_t$: $n_v = n_t/2 + 2 \leq n_t$, since $n_t$ is at least four.
Thus both the number of triangles and vertices on the boundary of the polygon are bounded by $n_* = \pi \gamma_*^2$. 
\end{proof}

On the class of bounded aspect ratio polyhedra with triangular facets of bounded aspect ratio, the standard interpolation error estimate holds.

\begin{theorem}
\label{th:boundedaspectratio3d}
There exists a constant $C > 0$ depending only upon $\gamma_*$ such that for any polyhedron $P\in\cP$ and all functions $u\in H^2(P)$,
\begin{equation}\label{eq:polybound3d}
\hpsn{u-I_P u}{1}{P} \leq C\, \diam(P) \hpsn{u}{2}{P}.
\end{equation}
\end{theorem}

\begin{proof}
The proof consists of two main ingredients. 
First we observe that the majority of the proof of Theorem~\ref{th:triangulationbasedbound} generalizes to this setting, leading to the conclusion that the interpolation error in the harmonic interpolant is no worse than that of the piecewise linear interpolation using any tetrahedralization that respects the triangle-faceted boundary.
Second, we show that a bounded aspect ratio tetrahedralization of any polygon in $\cP$ can formed using a single additional interior vertex at the incenter.
This is the same argument used in for the polygonal case in \cite[Theorem 2]{GRB2011}.

As usual, we restrict our analysis to polyhedra with unit diameter, without loss of generality. 
Following the proof of Theorem~\ref{th:triangulationbasedbound}, applying the Bramble-Hilbert lemma and Dirichlet's principle, we first estimate the interpolation error by
\begin{align}
\hpsn{u - I_P u}{1}{P}^2 & \leq \hpsn{u - p_u}{1}{P}^2 + \hpsn{p_u - I_P u}{1}{P}^2 \notag \\
 & \leq \hpsn{u - p_u}{1}{P}^2 + \sum_{T\in\cT}\hpsn{I_T\left(p_u -  u\right)}{1}{T}^2 \notag \\
 & \leq 2\hpsn{u - p_u}{1}{P}^2 + \sum_{T\in\cT}\hpsn{u - I_Tu}{1}{T}^2. \label{eq:poly3dint}
\end{align}
As before, $p_u$ is the linear polynomial from the Bramble-Hilbert lemma, but $\cT$ can be any tetrahedralization of $P$, possibly including additional interior vertices.
We require that each face of $P$ must be the face of one of the tetrahedra in $\cT$, i.e.\ there can be no mesh refinement along the boundary.

It remains to establish an estimate on $\hpsn{u - I_Tu}{1}{T}$. 
As the standard interpolation error for the linear interpolant on tetrahedra, this term can be estimated depending only upon the aspect ratio of $T$. 
Select $\cT$ to be the tetrahedralization formed by adding an additional vertex at an incenter of $P$ and forming tetrahedra by connecting this incenter to each of the faces. 
For each tetrahedron $T\in\cT$, a rigid body transformation will place the face of $T$ forming part of the boundary of $P$ in the $xy$-plane. 

By convexity, the aspect ratio bound on $\cP$ and the definition of the incenter, the distance from the incenter to the $xy$-plane (i.e., $h_*$ in Lemma~\ref{lm:tetqual}) is at least $1/\gamma_*$. 
Since $\cP$ requires the inradius of each face to be at least $1/\gamma_*$, Lemma~\ref{lm:tetqual} asserts the existence of a quality bound $\gamma_T$ depending only upon the original aspect ratio bound $\gamma_*$. 

The standard interpolation error estimate on tetrahedra asserts the existence of a constant $C_T$ depending only upon the aspect ratio of $T$ such that $\hpsn{u - I_Tu}{1}{T}^2 \leq C_T \diam(T) \hpsn{u}{2}{T}^2$.
Applying this to (\ref{eq:poly3dint}) and using the Bramble-Hilbert estimate completes the result:
\begin{align*}
\hpsn{u - I_P u}{1}{P}^2 
 & \leq 2\hpsn{u}{2}{P}^2 + \sum_{T\in\cT}C_T \hpsn{u}{2}{T}^2 = \left(2 + C_T\right)\hpsn{u}{2}{P}^2.
\end{align*}
\end{proof}

\subsection{Convex Polyhedra Without Bounded Aspect Ratio}
\label{ss:convex3d}

\begin{figure}[t]
     \begin{center}
\includegraphics[height=.2\textwidth]{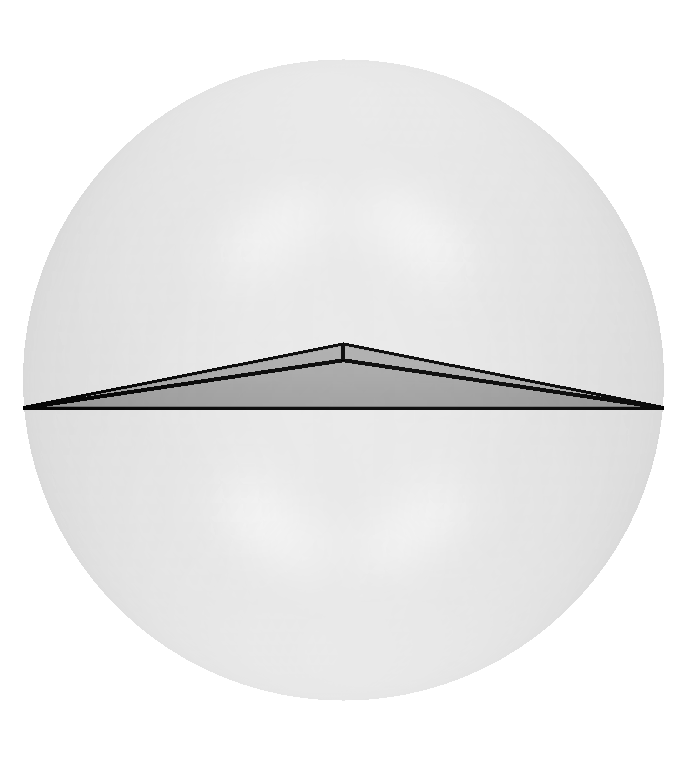}
\hspace{0.04\textwidth}
\includegraphics[height=.2\textwidth]{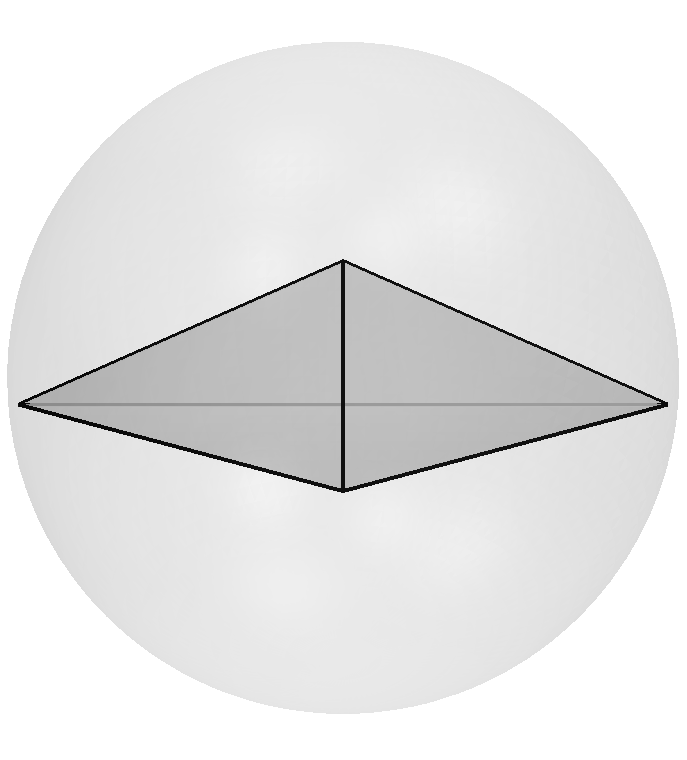}
\hspace{0.04\textwidth}
\includegraphics[height=.2\textwidth]{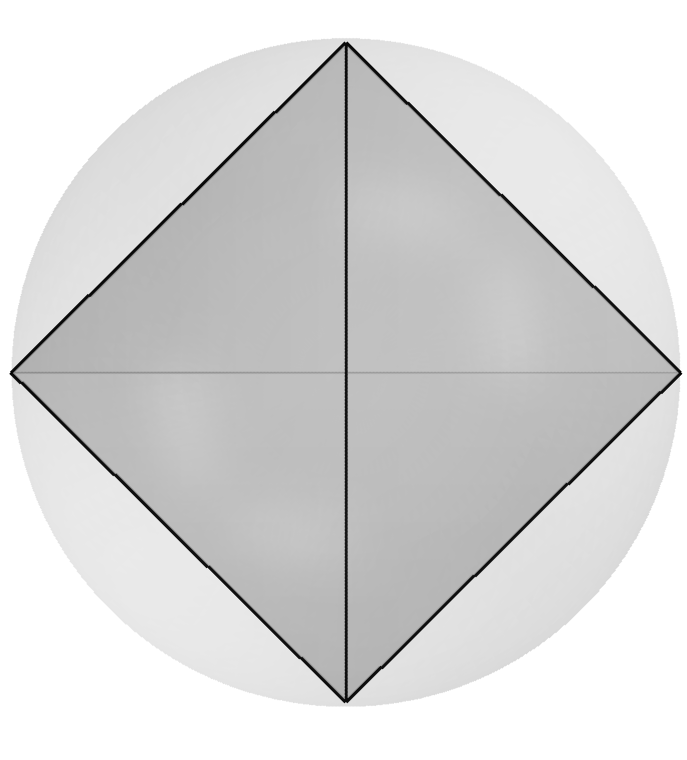}
     \end{center}
\caption{A sliver tetrahedron depicted from three different angles with its (modest sized) circumsphere.
}\label{fg:sliver}
\end{figure}

When translating the aspect-ratio independent error estimates on convex polygons in Sections~\ref{ss:esttri} and \ref{ss:convex} to three dimensions, much of the analysis is identical albeit without the connections to constrained Delaunay tetrahedralization. 
The inequalities (\ref{eq:poly3dint}) hold under a very limited set of assumptions, meaning that, in general, harmonic coordinates always perform at least as well as using a piecewise linear interpolant over a tetrahedralization of $P$.
Constructing an analog to Theorem~\ref{th:triangulationbasedbound} is not so straightforward, since interpolation error on tetrahedra is not bounded by circumradius, due to a class of tetrahedra called slivers.
Slivers are nearly coplanar tetrahedra with no short edges, all four vertices near a common circle, and circumradii proportional to the edge lengths.
These tetrahedra have poor interpolation properties, since they allow for disjoint edges that pass arbitrarily close to each other; see Figure~\ref{fg:sliver}. 
Sliver tetrahedra are not eliminated by the Delaunay construction and are a well known detriment to 3D Delaunay meshes; a substantial body of mesh generation literature is devoted to sliver removal, e.g.~\cite{CDEFT00,ELMSTTW00,EG02,Li03}.

Despite a great deal of interpolation error estimates for tetrahedra without a bound on aspect ratio (e.g.~\cite{Ja76,Kr92,Sh94}) the general error metric, called ``coplanarity'' in \cite{Ra12}, does not correspond to a natural geometric construction.
While Theorem~\ref{th:triangulationbasedbound} can be established in terms of a coplanarity measure, there is no apparent tool for generating tetrahedralizations of a given polyhedron with minimal coplanarity, even if it is assumed that all faces are triangles of high geometric quality.
Worse still, some polyhedra cannot be tetrahedralized without the introduction of additional vertices~\cite{Le1911,Sc28} including some convex polyhedra~\cite{Ra03}.
Accordingly, the constrained Delaunay tetrahedralization is not well-defined in general and typical approaches to ensure its existence involve adding additional vertices on the boundary~\cite{Sh08i}, which would interfere with the construction leading to estimate (\ref{eq:poly3dint}). 
While tetrahedralization does not provide a sharp characterization of interpolation error, many of the general principles still hold.
In particular, looking at the sharpness result in Section~\ref{ss:convex}, we now see that interpolation error is still guaranteed to grow whenever a vertex is very near the interior of a non-adjacent face.

Building the analog of Lemma~\ref{lm:convexcase}, consider a sequence $\left(P_i, \bv_i, f_i\right)_{i=1}^{\infty}$, where $P_i$ is a convex polyhedron, $\bv_i$ is a vertex of $P_i$ and $f_i$ is a face of $P_i$, under the following assumptions:
\begin{enumerate}[label=(B\arabic*)]
\item\label{as:diam13d} $\diam (P_i) = 1$;
\item\label{as:xyplane} $f_i$ lies in the $xy$-plane;
\item\label{as:zaxis} $\bv_i$ lies on the positive $z$-axis;
\item\label{as:nondegen3D} There exists a constant $c_v$ such that $\dist(\origin, \partial f_i) > c_v$ where $\partial f_i$ is the boundary of $f_i$ and $\origin$ is the origin;
\item\label{as:smalldist3D} $\dist(\bv_i, f_i) \rightarrow 0$ as $i\rightarrow \infty$. 
\end{enumerate}

\begin{figure}[ht]
     \begin{center}
\includegraphics[width=.4\textwidth]{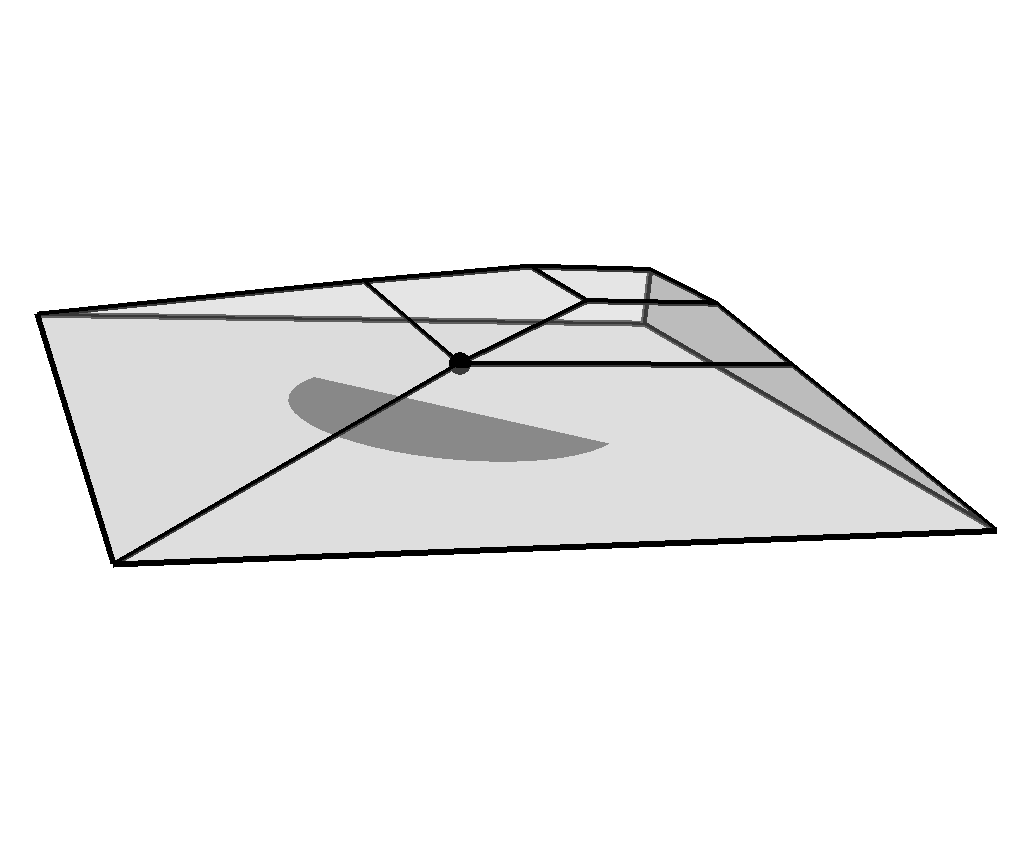}
\hspace{0.04\textwidth}
\includegraphics[width=.35\textwidth]{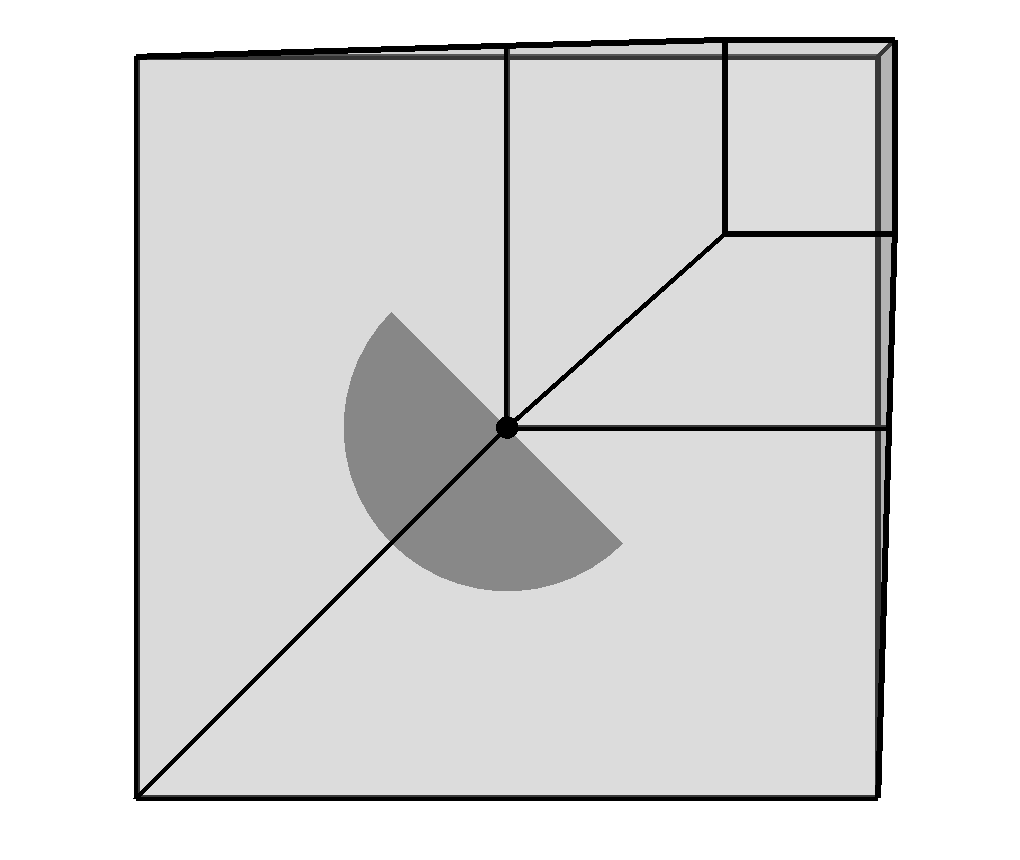}
     \end{center}
\caption{Example convex polyhedron containing a vertex (black dot) which is near the interior of a non-adjacent face. The integration domain where the gradient is known to be large is depicted with the dark gray half circle.
}\label{fg:convexThin3d}
\end{figure}

In this setting, the following lemma asserts that the standard error estimate cannot be established. 

\begin{lemma}\label{lm:convexcase3d}
For a sequence of polyhedron-vertex-face tuples $\left(P_i, \bv_i, f_i\right)_{i=1}^{\infty}$, satisfying (B1)-(B5),
\[
\lim_{i\rightarrow \infty} \frac{ \hpsn{u - I_{P_i} u }{1}{P_i} }{ \hpsn{u}{2}{P_i} } = \infty,
\]
for the function $u(x,y,z) = x^2+y^2$.
\end{lemma}
\begin{proof}
The proof follows the same steps as the 2D case, Lemma~\ref{lm:convexcase}. 
The isodiametric inequality again implies a uniform upper bound on the denominator, reducing the desired estimate to showing that the numerator grows without bound. 

Again, the integral in the numerator is restricted to the key portion near the origin, where large gradients occur in the $z$-direction:
\begin{align*}
\hpsn{u - I_{P_i} u }{1}{P_i}^2 & \geq \int_{P_i} \frac{\partial}{\partial z}\left(u - I_{P_i} u \right)^2 \notag \\
 & \geq \int_{D_i}^0 \int_0^{h_i(x,y)} \left(\frac{\partial}{\partial y}  I_{P_i} u(x,y)\right)^2 \dd z \dd A(x,y) .
\end{align*}
Above, $h_i(x,y)$ denotes the height of the boundary of the polygon above the $xy$-plane (which is well-defined near the origin above $f_i$). 
Here $D_i$ is the integration domain in the $xy-$plane: 
half of the ball centered at $(0,0)$ with radius $c_v^2/2$, where the plane supporting $P_i$ at $(0,0,h_i(0,0))$ is below $h_i(0,0)$; see Figure~\ref{fg:convexThin3d}. 
The inner integral is estimated using the piecewise linear nature of the interpolant on the boundary and we conclude that, in the integration domain,
\[
I_{P_i} u(x,y,h_i(x,y)) \leq c_v^2/\sqrt{2}.
\]
Then after applying Corollary~\ref{cor:minh1} to estimate the inner integral, we can estimate the remaining terms:
\begin{align*}
\hpsn{u - I_{P_i} u }{1}{P_i}^2 
& \geq \int_{D_i} \frac{\left(I_{P_i} u(x,y,h_i(x,y)) - I_{P_i} u(x,y,0) \right)^2}{h_i(x)}  \dd A(x,y)\\
 & \geq \int_{D_i} \frac{\left[ \left(1-1/\sqrt{2}\right) c_v^2 \right]^2}{\dist(\bv_i, f_i)} \dd A(x,y)\\
 & = \frac{\pi c_v^4}{8} \cdot \frac{\left[ \left(1-1/\sqrt{2}\right) c_v^2 \right]^2}{\dist(\bv_i, f_i)}.
\end{align*}
As $\dist(\bv_i, f_i)$ is the only remaining non-constant term, \ref{as:smalldist3D} ensures that the estimate grows without bound.
\end{proof}

\subsection{Non-convex Polyhedra}
\label{ss:nonconvex3d}

Non-convexity presents additional challenges when moving to three dimensions.
Remark~\ref{rm:sobo} suggested that vertices of non-convex polyhedra can approach non-adjacent faces without causing the interpolation error to blow up.
The lemma below formalizes the three dimensional analog of the example in Section~\ref{ss:nonconvexexample} demonstrating the a sharp characterization of polygons that admit bounded interpolation error estimates must include some shapes that are fundamentally different than those allowed in two dimensions.

\begin{lemma}
\label{lm:nonconvexcase3d}
There exists a sequence of polyhedron-vertex-face tuples $\left(P_i, \bv_i, f_i\right)_{i=1}^{\infty}$, satisfying \ref{as:diam13d}-\ref{as:smalldist3D} such that
\[
\lim_{i\rightarrow \infty} \frac{ \hpsn{u-I_{P_i} u}{1}{P_i} }{ \hpsn{u}{2}{P_i} } = C< \infty,
\]
for the function $u(x,y,z) = x^2 + y^2$ and some real number $C$.
\end{lemma}

\begin{figure}[ht]
\centering
  \begin{subfigure}[b]{0.4\textwidth}
    \includegraphics[width=\textwidth]{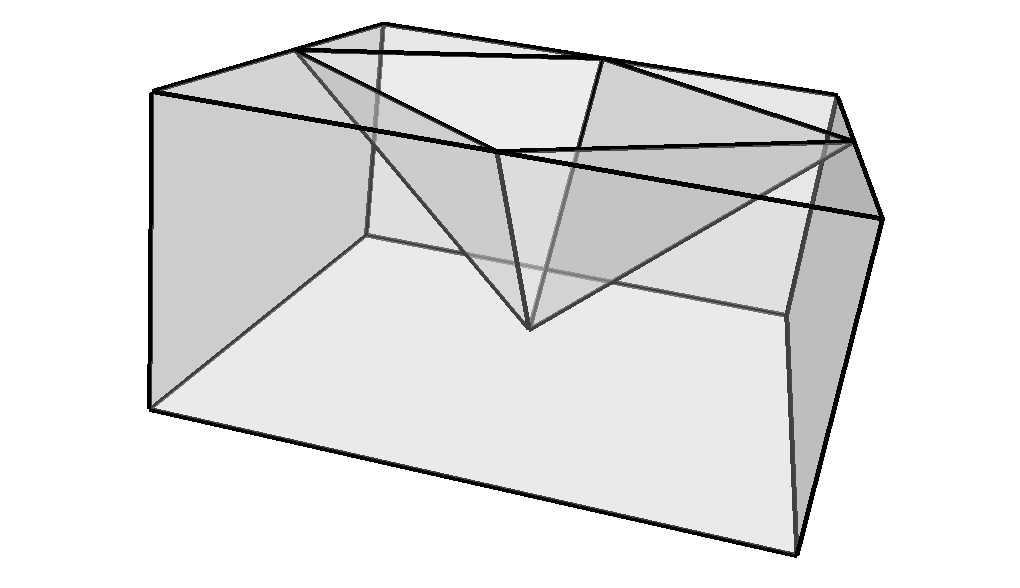}
    \caption{}
    \label{fg:nonconvex3d}
  \end{subfigure}
\hspace{0.08\textwidth}
  \begin{subfigure}[b]{0.22\textwidth}
    \includegraphics[width=\textwidth]{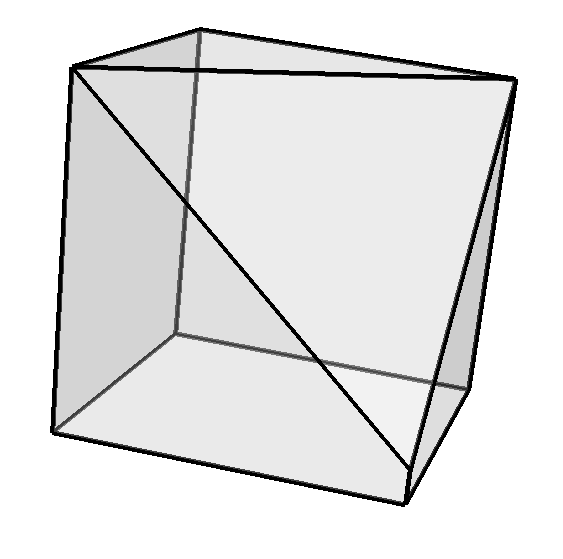}
    \caption{}
    \label{fg:convexPortion3d}
  \end{subfigure}
  \caption{
(\textsc{a}) A non-convex three-dimensional analog to the domains in Figure~\ref{fg:nonconvexPolygon} used in Lemma~\ref{lm:nonconvexcase3d}. The family of polyhedra involves allowing the vertex at the interior of the convex hull to approach the opposite face.  
(\textsc{b}) The convex quarter domain needed for analysis.
}
\label{fg:nonconvex3dfull}
\end{figure}

Most of the details are are a direct generalization of the construction given in Section~\ref{ss:nonconvexexample} to three dimensions. 
Here we outline the key points and differences from two dimensions.
The construction involves polyhedra shown in Figure~\ref{fg:nonconvex3d} with $13$ vertices: 
$(1,1,0)$, $(-1,1,0)$, $(1,-1,0)$, $(-1,-1,0)$, $(1,1,1)$, $(-1,1,1)$, $(1,-1,1)$, $(-1,-1,1)$, $(1,0,1)$, $(0,1,1)$, $(-1,0,1)$, $(0,-1,1)$ and $(0,0,\epsilon)$. 
The vertex $(0,0,\epsilon)$ is allowed to approach the non-adjacent face in the $xy$-plane. 

If all faces of these polyhedra are triangular, then $(0,0,\epsilon)$ is approaching a point on the boundary of a triangle and thus \ref{as:nondegen3D} is not satisfied. 
This can be corrected by adding vertices to the face of the polygon in the $xy$-plane so that $(0,0)$ lies in the interior of a triangle.
We ignore this minor detail, allowing polyhedra with non-triangular (in this case square) faces.

In the analysis, one fourth of the original domain is considered, namely, the portion of domain with $x > 0$ and $y > 0$ shown in Figure~\ref{fg:convexPortion3d}.
Symmetry of $u$ across the $xz$- and $yz$-planes ensures that analysis on this subdomain is sufficient: any integrals on the subdomain will be one fourth of the integral over the entire domain.
As in the family $\{Q_\eps\}$ in Section~\ref{ss:nonconvexexample}, the resulting domains are convex and have a similar ``uniform Lipschitzness.''

Rather than elaborate upon the details, we emphasize the key properties of the limiting case that are driving the result. 
The limiting boundary conditions have a discontinuity at a single point (the origin), but, in two dimensions, discontinuous functions are admitted in $H^{1/2}$, including these specific boundary conditions.
Thus, the limiting interpolant $I_{P_0}u$ is in $H^1$ and is in essence providing a uniform bound on the $H^1$-norm of interpolants for the cases with $\epsilon$ sufficiently near $0$.

The conclusion to be drawn from Lemma~\ref{lm:nonconvexcase3d} is as follows: the class of non-convex polyhedra that admit bounded interpolation errors is fundamentally different than all of the other geometry classes considered in this work.
For simplices (in any dimension), polygons (convex or non-convex) and convex polyhedra, a sequences of domains involving a vertex approaching an opposite edge/face causes an unbounded interpolation error while Lemma~\ref{lm:nonconvexcase3d} demonstrates that this is not the case for non-convex polyhedra.

\bibliographystyle{abbrv}
\bibliography{arxiv-hrmc}

\appendix

\section{Sobolev Spaces}\label{ap:sobolev}

For multi-index $\alpha = (\alpha_1, \alpha_2)$ and point $\bx = (x,y)$, define $\bx^\alpha := x^{\alpha_1} y^{\alpha_2}$, $\alpha ! := \alpha_1 \alpha_2$, $|\alpha| := \alpha_1 + \alpha_2$, and $D^\alpha u := \p^{|\alpha|} u/\p x^{\alpha_1}\p y^{\alpha_2}$.  
The Sobolev semi-norms and norms over functions defined on an open set $\Omega$ are defined by
\begin{align*}
\wmpsn{u}{m}{p}{\Omega}^p &:=  \int_\Omega \sum_{|\alpha| = m} |D^\alpha u(\bx)|^p \,{\rm d} \bx &{\rm and} & & \wmpn{u}{m}{p}{\Omega}^p &:= \sum_{0\leq k\leq m}\wmpsn{u}{m}{p}{\Omega}^p.
\end{align*}
In the case $p=2$, these norms define Hilbert spaces and a denoted $H^{m}(\Omega) := W^{m,p}(\Omega)$. 
The case $m=0$ reduces to space of integrable functions: the common $L^p$ spaces.

We now state a particular Sobolev embedding property that is most relevant to our analysis;
this result is phrased in much more generality in the literature on Sobolev spaces, e.g.,~\cite{Ad03,Le09}.

\begin{theorem}\label{th:morrey}(Morrey's inequality)
If $\Omega\in \R^n$ is a Lipschitz domain, and $m\geq n/p$, then all functions in $W^{m,p}(\Omega)$  are continuous. 
\end{theorem}

This theorem is valid for non-integral $m$, i.e.\ for fractional Sobolev spaces~\cite{DPV12}.
In our context, these fractional Sobolev spaces arise in the definition of trace spaces: the space of functions corresponding to the boundary values of the a Sobolev function. 
The trace theorem below asserts that such functions are well-defined in the simplest setting.

\begin{theorem}\label{th:trace}(Trace Theorem) \cite{M1987,Co88,Di96}
If $\Omega\in \R^n$ is a Lipschitz domain, then the trace operator,
\[ \textnormal{Tr}: H^1(\Omega) \rightarrow H^{1/2}(\partial \Omega) \]
is a well-defined, bounded, linear operator with a bounded right-inverse, i.e. there exists a constant $C > 0$ such that for all $u\in H^1(\Omega)$ and for all $g\in H^{1/2}(\partial \Omega)$,
\begin{align*}
\hpn{\textnormal{Tr } u}{1/2}{\partial \Omega} & \leq C \hpn{u}{1}{\Omega} \,\,\, \rm{and}\\
\hpn{\textnormal{Tr}^{-1} g}{1}{\Omega} & \leq C \hpn{g}{1/2}{\partial \Omega}.
\end{align*}
\end{theorem}

The typical theory of Sobolev spaces on Lipschitz domains involves covering the domain by patches and flattening the boundary on each patch to transform the domain locally into a half-space.
The following lemma is a technical detail used in Sections~\ref{ss:nonconvexexample} and \ref{ss:nonconvex3d} when arguing that certain domain-dependent constants are in fact domain independent for a specific family of domains.

\begin{lemma}\label{lm:flattening}
For any $\epsilon \in (0,1/4)$, the domain $Q_\epsilon$ (defined in Figure~\ref{fig:P-and-Q}) can be locally flattened using decomposition with a uniformly bounded Lipschitz constant that does not depend on $\epsilon$.
\end{lemma}

\begin{figure}[ht]
\[
\includegraphics[width=.3\textwidth]{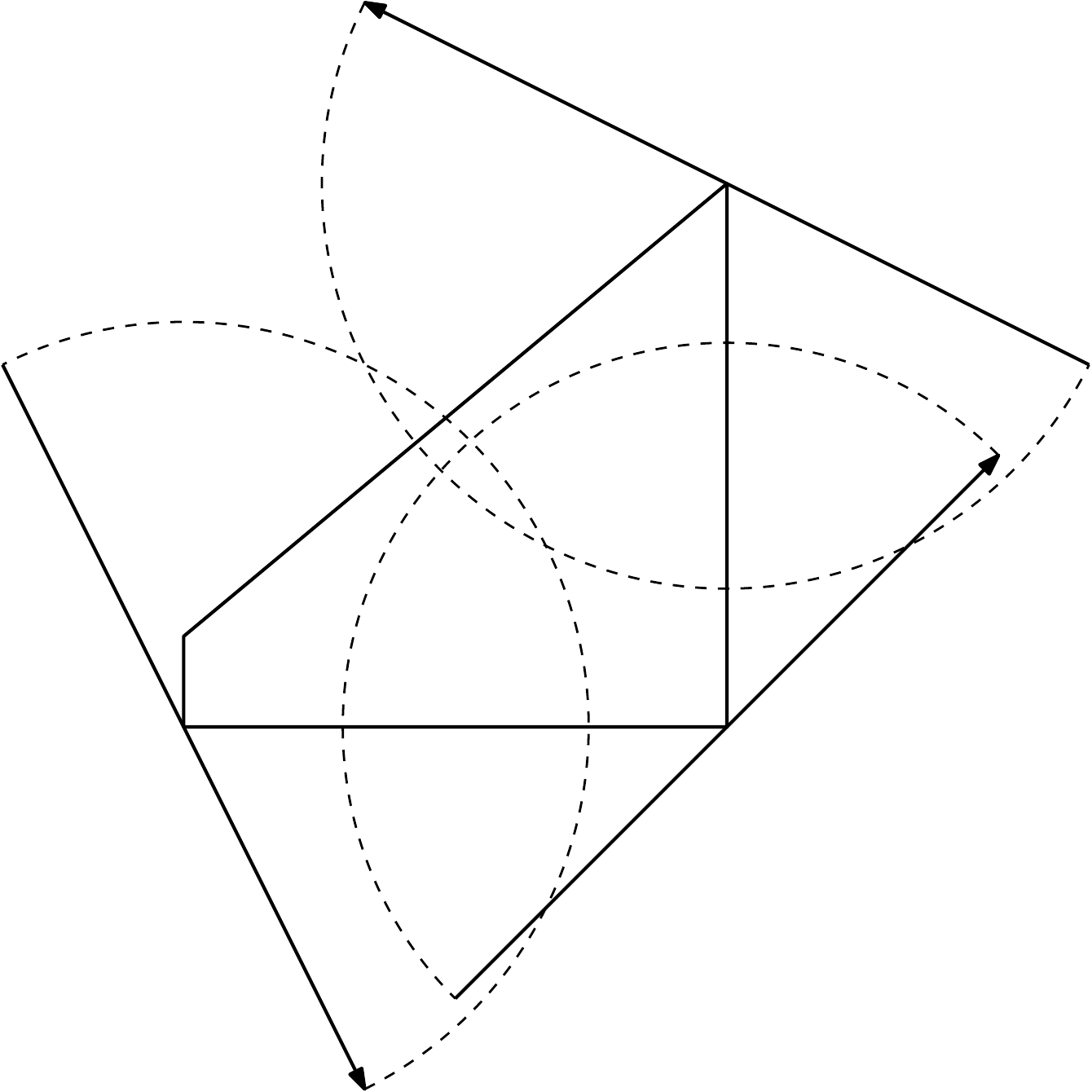}
\]
\caption{Local axes used for flattening the boundary of the domain $Q_\eps$ (see Figure~\ref{fig:P-and-Q} for context).  These local axes and patches can be used for all $\epsilon<1/4$.}\label{fg:domainFlattening}
\end{figure}

\begin{proof}
Consider flattening the domain in three different patches shown in Figure~\ref{fg:domainFlattening}. 
The origin of each of the local coordinates systems is located at one of the vertices of $Q_\epsilon$ that does not depend on $\epsilon$. 
We will show that the Lipschitz transformation for the patch centered at the origin is uniformly bounded. 
The patch centered at $(1,1)$ is similar and the patch centered at $(1,0)$ is simpler because locally the domain $Q_\epsilon$ is independent of $\epsilon$.

The local coordinate system is defined by the following (orthogonal) transformation:
\begin{align*}
\hx & = \frac{1}{\sqrt{5}}x - \frac{2}{\sqrt{5}}y;\\
\hy & = \frac{2}{\sqrt{5}}x + \frac{1}{\sqrt{5}}y.
\end{align*}
In this local coordinate system and for a ball of at radius at least $1/2$, the location of boundary of $Q_\epsilon$ can be written as a function $b$ of $\hx$:
\[
b(\hx) = \begin{cases} 
      2\hx & \textnormal{if } \hx> 0; \\
      -\frac{1}{2}\hx & \textnormal{if } -\frac{2}{\sqrt{5}}\epsilon \leq \hx\leq 0; \\
      \frac{3-\epsilon}{2\epsilon - 1}\hx + \frac{\sqrt{5}\epsilon}{2\epsilon - 1} & \textnormal{if } \hx < -\frac{2}{\sqrt{5}}\epsilon.
   \end{cases}
\]
For $\epsilon < 1/4$, we have $\vsn{b'(\hat x) } < 11/2$, so the derivative of $b$ is uniformly bounded with respect to $\eps$.
\end{proof}

\section{Dirichlet's Principle}
\label{ap:harmonic}

Dirichlet's principle asserts that harmonic functions have minimal $H^1$-norm among all functions with the requisite boundary conditions.
A precise statement and proof in a general setting can be found in Evans~\cite[Chapter 2, Theorem 17]{Ev98}. 
For reference, we state Dirichlet's principle in the notation of this paper.
\begin{theorem}[Dirichlet's principle]
\label{thm:dir-prin}
Assume $I_Pu\in C^2(P)$ solves (\ref{eq:interpde}) and let
\[\cA:=~\{w\in C^2(P)~:~w=g_u~\text{on $\p P$}\}.\]
Then 
\begin{equation}
\label{eq:dir-prin}\hpsn{I_P u}{1}{P} = \min_{w\in\cA}\hpsn{w}{1}{P}.
\end{equation}
Conversely, if $I_P u\in\cA$ satisfies (\ref{eq:dir-prin}), then $I_P u$ solves (\ref{eq:interpde}).
\end{theorem}

In 1D, the harmonic interpolant is the line fitting the values at the ends of an interval.
In this setting, Dirichlet's principle can be given as a more concrete inequality.
\begin{cor}
\label{cor:minh1}
Let $c_a, c_b \in \R$ and let $f_* : [a,b] \rightarrow \R$ be the linear function with $f_*(a) = c_a$ and $f_*(b) = c_b$. 
Then for all $f\in H^1([a,b])$ such that $f(a) = c_a$ and $f(b) = c_b$,
\[
\int_a^b \left(\frac{\dd }{\dd x} f(x)\right)^2 \dd x \geq \int_a^b \left(\frac{\dd }{\dd x} f_*(x)\right)^2 \dd x = \frac{\left(c_b - c_a \right)^2}{b-a}.
\]
\end{cor}

\section{Bramble-Hilbert Lemma}
\label{ap:bh}

The Bramble-Hilbert lemma~\cite{BH70} is the key estimate approximating Sobolev functions by polynomials in the appropriate context for the finite element method. 
While the basic finite element analysis only requires the theorem to hold on a specific reference element, the theorem can be established uniformly over convex domains~\cite{Ve99,DL04} which is essential to the analysis over general polygonal/polyhedral domains. 

\begin{theorem}\label{th:bh}
There exists a constant $C_{\ref{eq:bh}}$ depending only upon $n$ such that for any convex polytope of unit diameter $P\subset \R^n$ and any $u\in H^2(P)$, there exists an affine function $p_u$ such that
\begin{equation}\label{eq:bh}
\hpsn{u - p_u}{1}{P} \leq C_{\ref{eq:bh}} \hpsn{u}{2}{P}.
\end{equation}
\end{theorem}

\newpage
\section{Constrained Delaunay Triangulation}\label{ap:cdt}

Given a set of points $\cV = \{\bv_i\}$, the Delaunay triangulation consists of the set of triangles $\cT = \{T_j\}$ formed from vertices of $\cV$ satisfying the empty circumball property: 
the circumcircle of $T_j$ contains no points of $\cV$ in its interior. 
This uniquely defines a triangulation when the points of $\cV$ lie in general position, i.e., no four points share a single circle, but can be defined for any point sets using some arbitrary tie-breaking rules.

\begin{figure}[t]
\[\includegraphics[width=.3\textwidth]{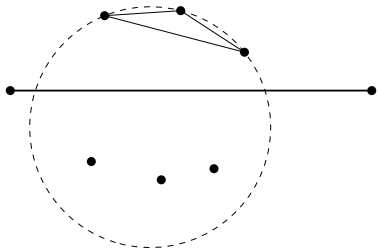}\]
\caption{A set of points and constrained segment for constrained Delaunay triangulation. The triangle shown is not part of the Delaunay triangulation of the points but is in the constrained Delaunay triangulation, since all vertices inside its circumcircle are not visible through the constrained segment.}\label{fg:visibility}
\end{figure}

Given a set of non-intersecting segments $\cS = \{S_k\}$ where each segment has end points in $\cV$, the constrained Delaunay triangulation~\cite{Ch89,Sh08i} is defined with a modified empty circumball criteria to create a triangulation with some similar properties while ensuring that the segments of $\cS$ appear in the triangulation.
Specifically, the empty circumball property is relaxed by allowing points in the triangle's circumball if they are not ``visible'' from the triangle.
A point is visible to a triangle if every line between that point and any point in the triangle does not cross a required segment; see Figure~\ref{fg:visibility}.

\begin{figure}[t]
\[\includegraphics[width=.5\textwidth]{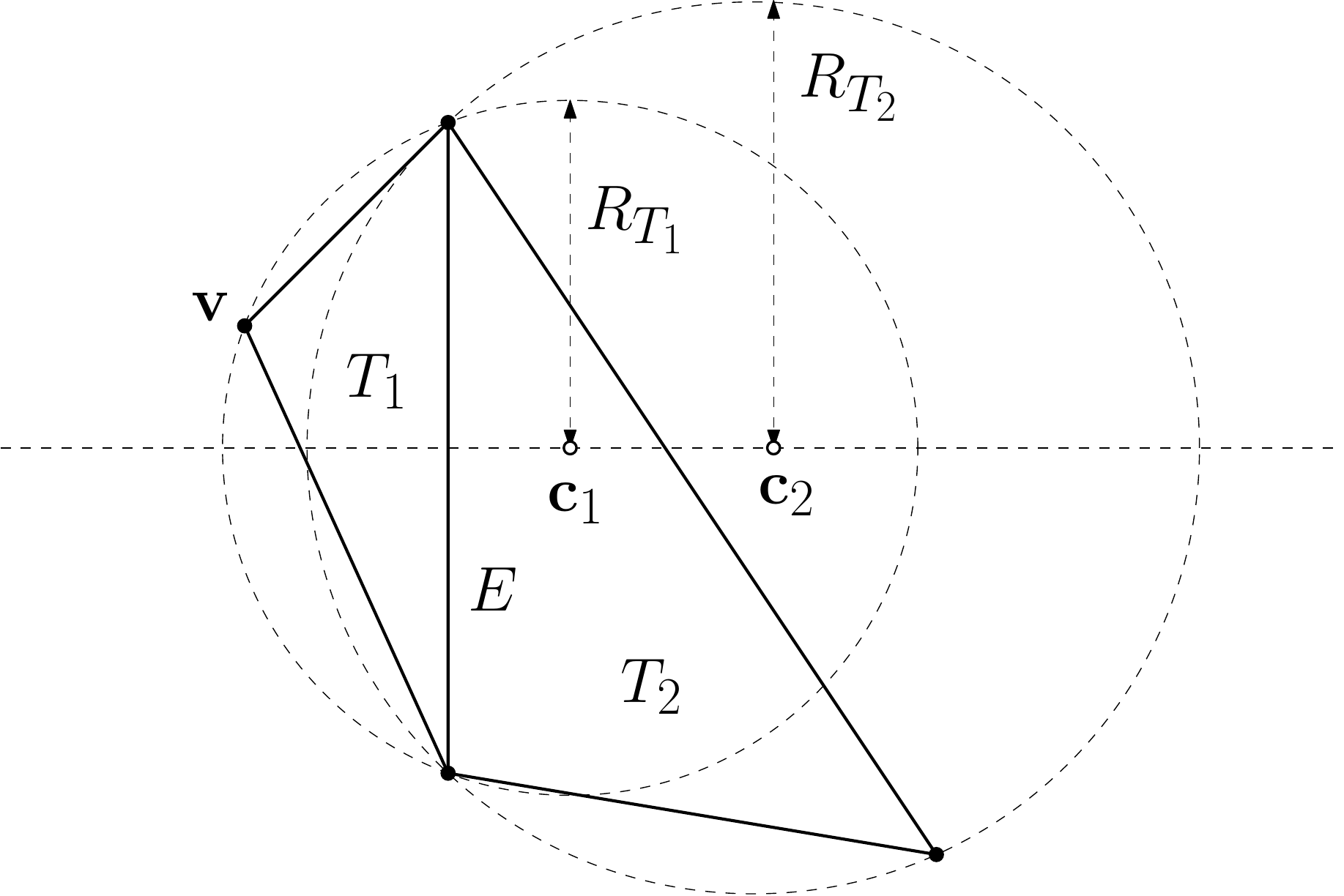}\]
\caption{In a Delaunay triangulation, the triangle adjacent to an obtuse triangle along the (unconstrained) longest edge always has a larger circumradius.}\label{fg:delaunayProof}
\end{figure}

\begin{lemma}\label{lm:delaunayProp}
Let $T_1$ and $T_2$ be adjacent triangles in a constrained Delaunay triangulation along a non-constrained edge $E$. If $T_1 $ is obtuse and adjacent to $T_2$ along edge $E$ then the circumradius of $T_1$ is smaller than the circumradius of $T_2$. 
\end{lemma}
\begin{proof}
Let $\bv$ denote the vertex of $T_1$ opposite edge $E$ and let $\bc_1$ and $\bc_2$ denote the circumcenters of $T_1$ and $T_2$, respectively; see Figure~\ref{fg:delaunayProof}. 
Since $T_1$ is obtuse, $\bc_1$ lies outside of $T_1$.
Since the circumcenters are equidistant from the endpoints of $E$, $\bc_1$ and $\bc_2$ lie on the line perpendicular to $E$.
If $\bc_2$ lies on the side of this line nearer to $T_1$, then the circumcircle of $T_2$ contains $\bv$ and thus is not in the constrained Delaunay triangulation.
So $\bc_2$ lies on the side of the line away from $T_1$ and thus the circumradius of $T_2$ is larger than that of $T_1$. 
\end{proof}

\section{Tetrahedron Quality}
\label{ap:tetqual}

Fixing constants $r_* > 0$ and $h_* > 0$, let $\cT_2$ denote the set of triangles in the $xy$-plane with inradius larger than $r_*$ and diameter at most $1$.  
Also define a set of points
\[
\cV = \left\{ (x,y,z) : \, \vsn{x^2 + y^2} \leq 1,\, 0 < h_* \leq z \leq 1 \right\}.
\]
Let $\cT_3$ denote the set of all possible tetrahedra with one face a triangle from $\cT_2$ and the fourth vertex from $\cV$. 
\begin{lemma}\label{lm:tetqual}
There exists a constant $\gamma_T < \infty$ depending only upon $r_*$ and $h_*$ such that the aspect ratio of every tetrahedron in $\cT_3$ is smaller than $\gamma_T$. 
\end{lemma}
A sketch of the proof of this lemma is as follows.
Since $r_* > 0$ and $h_*>0$, the aspect ratio of every tetrahedron in $\cV$ is strictly larger than $0$. 
Consider a sequence of tetrahedra $(T_i)$ such that limit of the aspect ratio approaches the infimum over $\cV$. 
Via compactness, we can extract a subsequence of this sequence such that the vertices converge. 
Thus by our original statement (all elements of $\cV$ have strictly positive aspect ratio), the lower bound on the aspect ratio of all tetrahedra in $\cV$ is strictly positive.

\end{document}